\documentclass{elsarticle}

\usepackage{amsmath,amssymb,amsthm,enumerate,multirow,subfigure,url}
\usepackage[normalem]{ulem}
\usepackage[usenames]{color}
\usepackage[utf8]{inputenc}
\usepackage[english]{babel}

\usepackage{url}

\usepackage{breakurl}
\usepackage[breaklinks]{hyperref}

\usepackage{here}
\usepackage{amssymb}
\usepackage[capitalize]{cleveref}
\usepackage{tikz}\usetikzlibrary{calc}\usetikzlibrary{positioning}
\usepackage{subfigure}
\usepackage[hmargin=2.5cm,vmargin=3cm]{geometry}
\usepackage[color=green!30]{todonotes}
\usepackage{comment}

\usepackage{perpage} 
\MakePerPage{footnote} 

\newtheorem{theorem}{Theorem}
\newtheorem{Question}{Question}

\newtheorem{lemma}[theorem]{Lemma}
\newtheorem{observation}[theorem]{Observation}
\newtheorem{corollary}[theorem]{Corollary}

\newtheorem{definition}[theorem]{Definition}

\newtheorem{conjecture}[theorem]{Conjecture}

\newcommand{\N}{\mathbb{N}}
\newcommand{\R}{\mathbb{R}}
\newcommand{\Z}{\mathbb{Z}}

\usepackage{algorithm,algpseudocode}

\begin{document}

\begin{frontmatter}

\title{The Generalized Double Pouring Problem: Analysis, Bounds and Algorithms}

\author[Umea]{Gerold J\"ager}
\author[Turku]{Tuomo Lehtilä}
\address[Umea]{Department of Mathematics and Mathematical Statistics,
University of Ume{\aa},\\
SE-901-87 Ume{\aa}, Sweden,\\
{\rm gerold.jager@umu.se}
}\address[Turku]{Department of Mathematics and Statistics,
University of Turku,\\
FI-20014 Turku, Finland,\\
{\rm tualeh@utu.fi}
}

\begin{abstract}
We consider a logical puzzle which we call double pouring problem,
which was original defined for $k=3$ vessels.
We generalize this definition to $k \ge 2 $ as follows.
Each of the $k$ vessels contains an integer
amount of water, called its value, where the values are $a_i$ for 
$i=1,2,\dots,k$ and the sum of values is $n$. 
A pouring step means pouring water
from one vessel with value $a_i$ to another vessel with value $a_j$, 
where $ 1 \le i \not= j \le k $ and $a_i \le a_j $.
After this pouring step the first vessel has value $2a_i$
and the second one value $a_j-a_i$.
Now the pouring problem is to find as few pourings steps
as possible to empty at least one vessel, or to show
that such an emptying is not possible (which is possible only in the case $k=2$).

For $k=2$ each pouring step is unique. We give a necessary and sufficient
condition, when for a given $ (a_1,a_2)$ with $a_1+a_2=n$ the pouring problem is solvable.
For $k=3$ we improve the upper bound of the pouring problem for 
some special cases. For $k \ge 4 $ we extend the known lower bound for $k=3$ 
and improve the known upper bound $\mathcal{O}((\log n)^2)$ for $k=3$ to $\mathcal{O}(\log n\log\log n)$.
Finally, for $k \ge 3$, we investigate values and bounds 
for some functions related to the pouring problem.

\end{abstract}

\end{frontmatter}

\section{Introduction}
\label{sec_intro}

\subsection{Water jugs pouring problem} 

A famous water pouring problem is the \emph{water jugs problem}. Its original version is formulated as follows.

\begin{itemize}
\item[] 
\emph{Given $2$ vessels with capacities $3$ and $5$, called \emph{value}, and an infinite water 
supply, can you precisely measure value $4$ by pouring steps from one vessel to 
the other one?}
\end{itemize}

This problem can easily be generalized by replacing $3$, $4$, $5$ by arbitrary numbers $a$, $b$ 
and $c$ and also to use $ k \ge 2 $ vessels $ a_1,a_2,\dots,a_k$ instead of $2$ vessels.  
The answer to this generalized problem is that this is possible if and only if 
$ \gcd (a_1,a_2,\dots,a_k) \, | \, c $~\cite{BSV02}. Another variant is to 
compute the minimum number of pouring steps for the original problem or one of 
its generalizations. In~\cite{ST08} it is shown that to find this minimum number 
for the generalized variant is $\mathcal{NP}$-hard. 
A heuristical solution is presented in~\cite{LZG05}.

\subsection{Sharing wine jugs problem} 

Another well-known water pouring problem is the \emph{sharing wine jugs problem}. 
Its original version is formulated as follows.

\begin{itemize}
\item[] 
\emph{Given $2$ vessels with capacities $a$ and $b$, called \emph{value}, and an 
infinite wine 
supply, can you divide the wine into equal parts by pouring steps from one vessel to 
the other one?}
\end{itemize}

Also this problem can easily be generalized to $ k \ge 3 $ vessels.
It was investigated for example in~\cite{AP76,HK23,MA94}.

\subsection{Double pouring problem} 

The problem of this work is similar to the both aforementioned pouring problems
and we call it \emph{double pouring problem}.
It was introduced as one of the 1971 All-Soviet-Union mathematical competition 
exercises as follows (see \cite{ASU71}, Exercise 148, 1971).

\begin{itemize}
\item[] 
\emph{The volumes of the water contained in each of $3$ big enough vessels are 
positive integers. You are allowed only to pour from one vessel to 
another the same volume of the water, that the destination vessel already contains.
Prove that you are able to empty one of the vessels.} 
\end{itemize}

It occurred in further mathematical competitions using different formulations
in USA and Canada (see \cite{WLP93}, Exercise B-6, 1993), 
in an IBM challenge (see \cite{IBM15}, 2015) and
in Germany (see \cite{BUI19}, Aufgabe 2, 2019).

\subsection{Previous results} 

Winkler presented in his book about mathematical puzzles \cite{Win04} a proof for the original exercise 
formulated above and also analyzed how many pouring steps have to be made at most
to ensure emptying one of the vessels. He provided a simple algorithm
with $ \mathcal{O}(n^2) $ pouring steps and an algorithm from Svante Janson
with $ \mathcal{O}(n \log n  ) $ pouring steps. 
Frei et al.\ investigated this problem 
in greater detail in~\cite{FRW21}. 
They improved the upper bound of \cite{Win04} from
$ \mathcal{O}(n \log(n) ) $ to $ \mathcal{O}((\log n )^2 ) $,
showed a lower bound $ \Omega(\log n  ) $ and presented
a class of hard instances, i.e., starting volumes for each of the three vessels.

\subsection{Generalized double pouring problem} 

In this work, we extend the double pouring problem from $k=3$ to 
$ k \ge 2 $ vessels, where each vessel has a volume and where each 
vector or values is called \emph{state}.
This means that the pouring steps are the defined similarly, but there
are $ \binom{k}{2} = k (k-1)/2$ possible pairs of vessels, where the pouring 
steps can take place.
The aim is again to empty one of the vessels, which is called \emph{being 
pourable}.

\subsection{Our results and organization of this work}

In Section~\ref{sec_def}, we give some notations and definitions which include 
introducing the generalized double pouring problem.

In Section~\ref{sec_two}, we consider the case $k=2$. 
We show that the problem essentially differs from the case 
$k\geq3$. Unlike in those cases, for $k=2$, there does not exist  for every initial state  a sequence 
of pouring steps which empties one of the vessels. We 
characterize all initial states which have a sequence of pouring steps leading 
to an empty vessel together with the exact number of pouring steps required in 
this sequence.

In Section~\ref{sec_three}, we consider the case $k=3$. 
We give two special cases of states for which we can improve 
the upper bound from $ \mathcal{O}((\log n )^2 ) $
to $ \mathcal{O}(\log n ) $.

In Section~\ref{sec_fixed}, we investigate values and bounds 
for some functions for the double pouring problem with at least three vessels.
In particular, we study the smallest and largest sum of values required for 
emptying one of the $k\geq3$ vessels when we need exactly $N$ pourings. 
For this purpose, we introduce four different functions dependent on $N$
and $k$. We present several monotony results in the parameters $N$ and $k$
and exact values for $N=1,2,3$ and a lower bound for $N=4$ for one of these 
functions.

In Section~\ref{sec_lower}, 
we show the same lower bound $ \Omega(\log n ) $ 
for the generalized double pouring problem, as in~\cite{FRW21}. 

In Section~\ref{sec_four}, 
as our main result, we improve the upper bound $ \mathcal{O}( (\log n )^2 ) $ 
from~\cite{FRW21}, which trivially holds also for $ k \ge 4 $,
to $ \mathcal{O}(\log n \log \log n ) $. 

In Section~\ref{sec_conc} we conclude and give suggestions for future work.

\section{Notations and Definitions}
\label{sec_def}

In this section, we introduce some basic notations and definitions.
Let $ \N$ be the set of positive integers and $ \N_0$ be the set of positive integers including $0$. We always write ``$ \log $'' for ``$\log_2 $''.
In the following, let $a,b,c,d,e,\dots$ be positive integers which stand for the initial content 
of the vessels of a pouring problem.

We begin by giving a definition introducing the pouring problem and a pouring 
step from one vessel to another.
\begin{definition}
\label{def_kvessel}
Let $ k \in \N $, and let $ A_1 ,A_2,\dots,A_k $ be a set of vessels
with value vector $ A= (a_1,a_2,\dots,a_k) \in \N_0^k $, called \emph{state}.
\begin{itemize}
\item[\bf (a)]  A \emph{pouring step} is a change from two vessels with entries $ a_i$ and $ a_j$
with $ 1 \le i \not= j \le n $ and $ a_i \le a_j $ to $ 2a_i$ and $ a_j -a_i$.
\item[\bf (b)] The \emph{pouring problem} is to find a sequence of pouring steps
so that at least one vessel has value $0$ in as few pouring steps as possible. 
\item[\bf (c)] Let $t\in \N_0$. \emph{Solving} the pouring problem in $t$ 
steps means emptying a vessel in at most $t$ pouring steps. 
\end{itemize}
\end{definition}

The following definition considers the number of pouring steps.
\begin{definition}
\label{def_pourable}
Let $A$ be a state.
\begin{itemize}
\item[\bf (a)]  For $p\in \N_0$, we say that $A$ is  \emph{$p$-pourable}, 
if a vessel can be emptied in at most $p$ pouring steps.
\item[\bf (b)]  For $p\in \N$, we say that $A$ is  \emph{exactly $p$-pourable}, 
if $A$ is $p$-pourable, but not $(p-1)$-pourable. 
\end{itemize}
\end{definition}

Some of our techniques rely on the notions of parity presented in the following definition.
\begin{definition}
Let $k\in\N_0$ and $a\in\Z \setminus\{0\}$. 
\begin{itemize} 
\item[\bf (a)]
We say that $a$ is $k$-{\rm even} if $2^k$ divides $a$. 
\item[\bf (b)]
We say that $a$ is $k$-{\rm odd} if $2^k$ does not divide $a$. 
\item[\bf (c)]
We say that $a$ is {\rm exactly} $k$-{\rm even} if it is $k$-even but not $(k+1)$-even.
In that case $a$ is said to have {\rm parity} $k$.
\end{itemize} 
\end{definition}

We will utilize the following parameters about the number of pouring steps of 
states.
\begin{definition}\label{defgh}
Let $N$ and  $k\geq3$ be positive integers. 
\begin{itemize}
\item[\bf (a)] 
Let $ g(N,k) $ be the smallest positive integer $n$ such that some state $A=(a_1,a_2,\dots,a_k)$ 
with $a_1+a_2+\cdots+a_k=n$ is not $(N-1)$-pourable. 
\item[\bf (b)] Let $ g'(N,k) $ be the smallest positive integer $n$ such that all states $A=(a_1,a_2,\dots,a_k)$ 
with $a_1+a_2+\cdots+a_k=n$ are $N$-pourable and at least one state is exactly $N$-pourable  (or equivalently, and at least one state is not $(N-1)$-pourable). 


\item[\bf (c)] 
Let $ h(N,k) $ be the largest positive integer $n$ such that all states $A=(a_1,a_2,\dots,a_k)$ 
with $a_1+a_2+\cdots+a_k=n$ are $N$-pourable.

\item[\bf (d)] 
Let $ h'(N,k) $ be the largest positive integer $n$ such that all states $A=(a_1,a_2,\dots,a_k)$ 
with $a_1+a_2+\cdots+a_k=n$ are $N$-pourable and at least one state is exactly $N$-pourable 
(or equivalently, and at least one state is not $(N-1)$-pourable).

\end{itemize}
\end{definition}

\section{Solvability for two vessels}
\label{sec_two}

In this section, we consider the pouring problem when we have $k=2$ vessels. As 
we observe later in this section, the case with only two vessels radically 
differs from the cases with at least three vessels. More concretely, we will see that for 
two vessels there exist some values of $n$ and initial states from which we 
cannot empty either of the two vessels. Due to this property, we have required 
that $k\geq3$ in Definition~\ref{defgh}.

\begin{definition}
\label{def_2vesselprob}
\begin{itemize}
\item[\bf (a)]
Define the {\rm two-dimensional pouring function}
$f: \N^2 \to \N_0^2 $ as follows
\begin{eqnarray*}
f(a,b) & := & \left\{ 
\begin{array}{ll}
(2a,b-a) & \mbox{if $ a\le b$,}
\\
(a-b,2b) & \mbox{otherwise},
\end{array}
\right.
\end{eqnarray*}
where $ a,b \in \N$.
\end{itemize}
\item[\bf (b)] A state $ (a,b) \in \N^2 $ is pourable,
if $ p \in \N $ exists such that
$ f^p(a,b) = (a+b,0)$. 
Then we say that the state $(a,b)$ is $p$-pourable.
\end{definition}

We start with a simple observation which we will apply in Section~\ref{sec_four} for four vessels.
\begin{observation}
\label{obs_bas}
If we pour once or several times between two vessels $A$, $B$ with sum of values $ n=|A| + |B| $, then the maximum 
of both values is at least $n/2$.
\end{observation}

In the following two lemmas, we present some pourable and some non-pourable initial states.

\begin{lemma}\label{Lem_2vesseleven}
Let $a,b\in \N$ and $a+b$ be odd. 
Then $ (a,b) $ is not pourable.
\end{lemma}

\begin{proof}
By applying the function $f$ recursively, the sum of the two values remains $a+b$.
If $ f^{k}(a,b) = (a+b,0) $, then   
$ f^{k-1}(a,b) = ((a+b)/2,(a+b)/2) $ follows. This is not possible,
as $a+b$ is odd.
\end{proof}

\begin{lemma}\label{Lem_2vesselgcd}
Let $a,b\in \N$. The state
$(a,b)$ is pourable 
if and only if $(a/\gcd(a,b),b/\gcd(a,b))$ is pourable.
\end{lemma}

\begin{proof}
Let  $ l := \gcd(a,b)$ and w.l.o.g., let $ a \le b $.
It holds:
\begin{eqnarray*}
f(a/l,b/l) \; = 
\; (2a/l,b/l-a/l)
\; = \; (2a, b-a )/l \; = \; f(a,b)/l. 
\end{eqnarray*}
The assertion follows.
\end{proof}


Now, we are ready to characterize all pourable states for $k=2$.

\begin{theorem}
\label{the_2vesselOnlyif}
Let $a,b\in \N$. The state
$(a,b)$ is pourable if and only if $(a+b)/\gcd(a,b)=2^k$ for some $k\in \N$
holds.
\end{theorem}

\begin{proof}
By Lemma~\ref{Lem_2vesselgcd}, we may assume that $\gcd(a,b)=1$. 
\begin{description}
\item[``$\Rightarrow$''] 
Let $(a,b)$ be pourable. Assume to the contrary 
that $a+b\neq 2^k$ for any $k\in \N$.

Thus, there must be an odd prime $p$ which divides $ a+b$.
If $p$ divided $a$, it would also divide $b$, a contradiction to 
$\gcd(a,b)=1$. 
The state $(a,b)$ can be pourable only if for some $ k \in \N $ with 
$ f^{k-1}(a,b) = (a',b') $ and
$ f^{k}(a,b) = (a'',b'') $ it holds that $ p $ divides neither $a'$
nor $b'$, but divides both $a''$ and $b''$.
This cannot hold, as one of $\{a'',b''\}$ is double of one of $ \{a',b'\} $,
and $p$ is an odd prime.

\item[``$\Leftarrow$''] 
Let $a+b= 2^k$ for some $k\in \N$.
Let us prove the claim with induction on $k$. 

\begin{description}
\item[$\mathbf{k=1:}$] 
The claim is clear. 
\item[$\mathbf{k \to k+1:}$] 
Assume that the claim holds for $k\in \N$. Let 
$a+b=2^{k+1}$ and w.l.o.g., let $a \le b$. Notice that $a$ and $b$ are both odd, 
as $ \gcd(a,b) = 1 $. 
We have $ f(a,b) = (2a,b-a)$. Now, both $2a$ and $b-a$ are even. Thus, Lemma 
\ref{Lem_2vesselgcd} implies that 
$(2a,b-a)$ is pourable if and only if $(a,\frac{b-a}{2})$ is pourable. Moreover, we 
have $a+\frac{b-a}{2}=2^k$ and hence, 
$(a,\frac{b-a}{2})$ is pourable by induction assumption.
Thus, also $(a,b)$ is pourable.\qedhere
\end{description}
\end{description}
\end{proof}

We end this section with a theorem, for $k=2$, giving the exact number of pourings 
whenever the initial state is pourable. We note that this result is significantly 
stronger than Theorem~\ref{thmBoundPourings3} for $k=3$ in which we do not know even the 
correct magnitude for the required number of pouring steps.

\begin{theorem}
\label{Lem_2vesselStepNumber}
Let $a,b,k\in \N$
and $(a+b)/\gcd(a,b)=2^k$. Then $(a,b)$ is pourable in exactly $k$ steps.
\end{theorem}
\begin{proof}
Because of $(a+b)/\gcd(a,b)=2^k$, the state $(a,b)$ is pourable by 
Theorem~\ref{the_2vesselOnlyif}. Let $t$ be the number of pouring steps.
Note that all pouring steps are 
uniquely defined. Let $f^{i}(a,b)=(a_i,b_i)$ for $ i=0,1,\dots,t$.
Thus, $(a,b)=(a_0,b_0) $ holds.
We prove the following claim.

\medskip

\noindent
{\bf Claim:}
$\gcd(a_i,b_i)=2^i \gcd(a,b)$  for $i=0,1,\dots,t$.

\medskip

\begin{description}
\item[\bf Proof (Claim):] Induction on $t$.

\begin{description}
\item[$\mathbf{t = 0:}$] 
This holds by assumption.

\item[$\mathbf{t-1 \to t:}$] 
Assume that the claim is true for $i \le t-1$, i.e., 
$\gcd(a_i,b_i)=2^i \gcd(a,b)$ holds. 
W.l.o.g.,  assume that $a_i \le b_i$ and let 
$a_i=\gcd(a_i,b_i)a_i'$ and $b_i=\gcd(a_i,b_i)b_i'$ for $i=0,1,\dots,t$. 
We have $a_{i+1}=2a_i$ and 
$b_{i+1}=b_i-a_i=\gcd(a_i,b_i)(b_i'-a_i')$. Observe that we have 
$a_i'+b_i'=a+b=2^k \gcd(a,b)$ and 
$\gcd(a_i',b_i')=1$. So both $a_i'$ and $b_i'$ are odd and $b_i'-a_i'$ is even. 
Thus, by induction assumption $\gcd(a_{i+1},b_{i+1})=2\gcd(a_i,b_i)=2^{i+1} 
\gcd(a,b) $ follows. 
\end{description}\end{description}

Let $i\in\{0,1,\dots,k\}$. 
By the claim, it holds that
$(a_i+b_i)/\gcd(a_i,b_i)
=(a+b) / (\gcd(a,b) \cdot 2^i) 
= 2^{k-i}$. Then the pouring process is finished
if and only if $a_t=a+b$ and $b_t=0$. This is equivalent
to $2^t \gcd(a,b)=\gcd(a_t,b_t) = a_t+b_t = a+b $ and equivalent to $ t=k$.
Thus, $(a,b)$ is pourable in exactly $k$ steps.
\end{proof}

\section{Upper bounds for special cases of three vessels}
\label{sec_three}

In this section, we first recall some known results for three vessels and then 
introduce some new results about the functions $ g(N,3) $ and 
$h(N,3)$ and about some initial states which need only few pouring steps.

\subsection*{Previous results}
We begin by presenting some algorithms for the pouring problem for three vessels. 
Algorithm \ref{alg_janson} originates from \cite[Pages 79, 84, 85]{Win04} and 
Algorithm \ref{alg_frei} from \cite{FRW21}. The two algorithms are utilized 
later in Section~\ref{sec_four}.

\begin{algorithm}[!th]
	\caption{Round of Janson's algorithm}
	\begin{algorithmic}[1]
	\Require  State $(a,b,c)$ with $a,b,c\in \N$ and $a\le b\le c$.
	\Ensure New state $(a',b',c')$ in which $b'<a$. 
\State Let $p= \lfloor b/a \rfloor $. 
\State Let $p=\sum_{i=0}^{\lfloor\log p\rfloor}p_i 2^i$ where each $p_i\in\{0,1\}$.
\For{$i=0,1,\dots, \lfloor\log p\rfloor$}
\If{$p_i=1$}
\State Pour from $B$ to $A$.
\Else{
Pour from $C$ to $A$.
}\EndIf
\EndFor
			\end{algorithmic}
	\label{alg_janson}
	\end{algorithm}

\begin{algorithm}[th!]
	\caption{Round of Frei's algorithm}
	\begin{algorithmic}[1]
	\Require  State $(a,b,c)$ with $a,b,c\in \N$ and $a\le b\le c$.
	\Ensure New state $(a',b',c')$ in which $a'\le a/2$ or $ b' < a/2 $. 
\State Let $p= \lfloor b/a \rfloor $, $ q = \lceil b/a \rceil $.
\State Let $r_1= b-pa $, $ r_2 = qa-b $.
\If{$r_1\le r_2$}
\State Apply Algorithm~\ref{alg_janson}.
\Else{
\State Let $q= \sum_{i=0}^{\lfloor\log (q)\rfloor}q_i 2^i$ where each $q_i\in\{0,1\}$.
\For{$i=0,1,\dots, \lfloor\log q\rfloor-1$.}
\If{$q_i=1$}
\State Pour from $B$ to $A$.
\Else{
Pour from $C$ to $A$.
}\EndIf
\EndFor
\State Pour from $A$ to $B$.
\EndIf
}
\end{algorithmic}
\label{alg_frei}
\end{algorithm}

Note that by applying Algorithm \ref{alg_frei} multiple times, Frei et al. 
have shown the following result.

\begin{theorem}[\cite{FRW21}, Theorem 1]
\label{thmBoundPourings3}
Let $(a,b,c)\in \N^3$ be a state and $n=a+b+c$. Then solving the 
pouring problem needs at most $ (\log n)^2$ pouring steps.
\end{theorem}

The following two lemmas are based on~\cite{FRW21}.

\begin{lemma}
\label{LemSvanteStateshift}
Let $a,b,c,q\in \N$, $a\le b\le c$, $b=qa+r$, 
where $0\le r<a$, 
and let $ h = \lfloor \log q \rfloor $.
Then the state $(a,b,c)$ can be transformed 
by pouring to the state $(2^ha,2^ha+r,c')$, for some $c' \in \N$, in $h $ steps.
\end{lemma}
\begin{proof}
Consider the binary representation of $q$, i.e., let 
$q = \sum_{i=0}^{h} q_i 2^i $, 
where $q_i \in \{0,1\} $ for $ i=0,1,\dots,h$, and $q_h=1$. 
We have $b=\left(\sum_{i=0,q_i=1}^h 2^{i}\right)a+r$. We do $h$ pouring steps, first for $q_0$, then for $q_1$ 
and so on, and finally for $q_{h-1}$. 
Let $i \in\{0,1,\dots, h-1\}$.
If $q_i=1$, then we 
pour from $B$ to $A$, and if $q_i=0$, then we pour from $C$ to $A$. 
This leads to the final state $(2^ha,2^ha+r,c')$ for some $c' \in  \Z$. Notice that $c'$ is positive 
since $\left(\sum_{i=0}^{h-1}2^{i}\right) a < qa\le b\le c$ and $ c' \ge c - 
\left(\sum_{i=0}^{h-1}2^{i}\right)a \ge 1 $.
\end{proof}

\begin{lemma}
\label{comp_jan_frei}
Both Algorithm~\ref{alg_janson} and Algorithm~\ref{alg_frei} use at most
$ \log\left( \frac{b}{a} \right) +2 $ pourings.	
\end{lemma}
\begin{proof}
Observe that in both algorithms we have $p= \lfloor b/a \rfloor $ and in 
Algorithm~\ref{alg_frei} we also have $ q	= \lceil b/a \rceil $. 
Because of $ p \le q $, in both algorithms we have at most 
$ \lfloor \log q \rfloor + 1  
\le \lfloor \log \left( \frac{b}{a}+1 \right)\rfloor+1  
\le \log \left( \frac{b}{a} \right) + 2 $ pourings.
Thus, the claim follows.
\end{proof}

\begin{theorem}
[\cite{FRW21}, Theorem 2]
\label{thek3hBound}
For each positive integer $n$, there exists an  initial state $A=(a,b,c)$ with 
$a+b+c=n$, which needs at least $\lceil\log((n+1)/5)\rceil$ pourings  
to solve the pouring problem.
\end{theorem}

\subsection*{Our results}

As a corollary of Theorem~\ref{thmBoundPourings3}, we can prove
a lower bound for $g(N,3)$.  

\begin{corollary}
\label{corgN3}
It holds that $g(N,3)\geq 2^{\sqrt{N}}$.
\end{corollary}

\begin{proof}
By Theorem~\ref{thmBoundPourings3}, we have  $N\leq (\log g(N,3))^2$. Hence, 
$g(N,3)\geq 2^{\sqrt{N}}$ holds.
\end{proof}

Similarly, as a corollary of Theorem~\ref{thek3hBound}, we can prove
an upper bound for $h(N,3)$.  

\begin{corollary}
\label{corhN3}
It holds that $h(N,3)\leq 5\cdot2^{N}-1$.
\end{corollary}

\begin{proof}
By Theorem~\ref{thek3hBound}, we have $N \geq \log((h(N,3)+1)/5)$. Hence, 
$h(N,3)\leq 5\cdot2^{N}-1$ holds.
\end{proof}

Based on computational results (see Table~\ref{tab:hnk} of 
Section~\ref{sec_fixed}), we present the following 
conjecture for $h(N,3)$. Note that Frei et al. \cite{FRW21} have shown that 
there exists infinitely many initial states satisfying the following conjecture. 
\begin{conjecture}\label{con:hn3}
For each positive integer $N$, we have $h(N,3)= 5\cdot2^{N-1}$.
\end{conjecture}

In the following theorem, we show that for some initial states, we need
only a small number of pourings.

\begin{theorem}
\label{TheSvanteAlg}
Let $a,b,c\in \N$ and $a$ be divisor of $b$, and let $n=a+b+c$. Then it holds:

\begin{itemize}
\item[\bf (a)]
Let $a\le b\le c$ and $b=qa+r$, where $ q \in \N $ and $0\le r<a$. 
Then the state $(a,b,c)$ is pourable in $(r+1) \lfloor \log n \rfloor $ steps. 
\item[\bf (b)]
The state $(a,b,c)$ is pourable in $ \lceil \log n \rceil $ steps. 
\end{itemize}
\end{theorem}

\begin{proof}
\begin{itemize}
\item[\bf (a)]
Let $q = \sum_{i=0}^{h} q_i 2^i $, 
where $q_i \in \{0,1\} $ for $ i=0,1,\dots,h$, and $q_h=1$. 
By Lemma~\ref{LemSvanteStateshift}, we can obtain the state $(2^ha,2^ha+r,c')$ 
for some $ c' \in \N $ in
$\lfloor \log q \rfloor $ steps from the state $(a,b,c)$. Furthermore, the state 
$(2^{h+1}a,r,c')$ can be 
obtained in at most $ 1+\lfloor \log q \rfloor \le 1+\lfloor \log b \rfloor 
\le \lfloor \log n \rfloor $ steps. 
Now we repeat this process.
As in each repetition the smallest value decreases, we need at most $r+1$ rounds each 
taking at most $\lfloor \log n \rfloor $ steps to arrive at a state with an empty vessel.
In total, this leads to $(r+1) \lfloor \log n \rfloor $ pouring steps. 
\item[\bf (b)]
This follows from (a) with $r=0$.\qedhere
\end{itemize}
\end{proof}

In the following lemma and theorem, we give a link between the state $(a,b,c)$ 
having a parity of $a+b+c$ that is $t$-even and the number of pouring steps 
needed. In particular, the larger the value $t$ is, the less pouring steps are
needed. Later 
in Section \ref{sec_four}, we use this idea to give an efficient algorithm for 
solving the pouring problem for four vessels.

\begin{lemma}
\label{lem_3vesselgcd}
Let $a,b,c\in \N$. 
Let the state $(a',b',c')$ be reachable from the state $(a,b,c)$ 
with some sequence of pourings. Then $\gcd(a',b',c')=2^t\gcd(a,b,c)$ for some 
$t\in \N_0$.

\end{lemma}
\begin{proof}
Let $ g:= \gcd(a,b,c) = 2^l \cdot q $, where $q$ is odd, 
and let $ g' := \gcd(a',b',c') $.
First we show that $g \, | \, g'$.

Clearly, $g$ divides $a$, $b$ and $c$.
W.l.o.g., assume that we pour from the second vessel to the first one. This leads to
the state $(2a,b-a,c)$. 
Here $g$ divides each of $2a,b-a$ and $c$. 
A similar argument holds after each further pouring step.
Thus, $g $ divides each of $ a',b'$ and $c'$,
and thus also their $\gcd$. So $g \, | \, g'$ follows. 

Next assume that $g':= 2^h \cdot q \cdot q' $ for $h\in\N$ with $ h \ge l $, and $q' 
\not= 1 $ is odd.
Then $q q'$ divides each of $a'$, $b'$ and $c'$ but does not divide each of $a$, $b$ and $c$. 
Consequently, (possibly after reordering) there exists a state
$(a_1,b_1,c_1)$ and a state after pouring $(2a_1,b_1-a_1,c_1)$ such that  $qq'$ does 
not divide both $a_1$ and $b_1$ but divides both $2a_1$ and $b_1-a_1$. However, this is 
not possible since $q$ and $q'$ and thus also $q \cdot q'$ are odd. Thus, $ 
q'=1$ and the assertion holds with $ t := h-l$.
\end{proof}

Note that the previous lemma can also be interpreted in the following way: 
The only factors of $\gcd(a,b,c)$ which may change when we pour from one vessel 
to another are the multiplicities of 2. Furthermore, the parity of $\gcd(a,b,c)$ 
may only remain unchanged or increase.

\begin{theorem}
\label{thm_3vesselStepNumber2n}
Let $a,b,c,l\in \N$ with $ n=a+b+c$ and $n/\gcd(a,b,c)=2^l$.
Then it holds:
\begin{itemize} 
\item[\bf (a)]
The state $(a,b,c)$ is pourable and after exactly 
$l$ steps we have two empty vessels.
\item[\bf (b)]
The state $(a,b,c)$ is pourable in at most $ \lfloor\log n\rfloor$ steps. 
\end{itemize} 
\end{theorem}

\begin{proof}
\begin{itemize} 
\item[\bf (a)]
We have $(a+b+c)/\gcd(a,b,c)=2^l$. Set $a_0=a/\gcd(a,b,c)$, 
$b_0=b/\gcd(a,b,c)$ and $c_0=c/\gcd(a,b,c)$. Notice that exactly two of $a_0$, 
$b_0$ and $c_0$ are odd since $a_0+b_0+c_0=2^l$ and $\gcd(a_0,b_0,c_0)=1$. 
W.l.o.g., $a_0\le b_0$ and $a_0$ and $b_0$ are odd. 
Hence, we can obtain the state $(2a_0,b_0-a_0,c_0)$ with a single pouring step. 
Moreover, $\gcd(2a_0,b_0-a_0,c_0)=2$, as each of $2a_0$, $b_0-a_0$ and 
$c_0$ is even, and $4$ does not divide $2a_0$ and by Lemma~\ref{lem_3vesselgcd}, 
no other positive integer larger than one can divide each of $2a_0$, $b_0-a_0$ and $c_0$. Thus, we may 
consider $a_1=a_0/2$, $b_1=b_0/2$ and $c_1=c_0/2$. 

We have $a_1+b_1+c_1=2^{l-1}$ and $\gcd(a_1,b_1,c_1)=1$. We continue in this way 
iteratively, each round pouring between the two odd vessels. Since the sum 
$a_i+b_i+c_i$ halves during each iteration, at some point, let us say on the $t$-th 
iteration, we reach the state $(a_t,b_t,c_t)$ where two of these values are 
identical. The last point this might occur is when $a_t+b_t+c_t=4$, 
where $a_t=b_t=1$ and $c_t=2$. 
After this, we have the state $(0,b_{t+1},c_{t+1})$ and $0 + b_{t+1}+c_{t+1}=2^{l-t-1}$. 
By Theorem~\ref{Lem_2vesselStepNumber}, $ (b_{t+1},c_{t+1})$ 
is pourable in $l-t-1$ steps. Thus, 
we obtain two empty vessels in $l$ steps.
\item[\bf (b)]
This follows from (a), as $ \log n = \log (a+b+c) = \log (\gcd(a,b,c)) + l$
and thus $ l \le \lfloor \log n \rfloor $.\qedhere
\end{itemize} 
\end{proof}

\section{Analysis for fixed number of pourings for at least three vessels}
\label{sec_fixed}

Whereas Frei et al.~\cite[Table 1]{FRW21} found the values of $g(N,3)$ for  $1\leq 
N\leq 14$, Tromp~\cite[A256001]{oeis25} has calculated those values for 
$15\leq N\leq 18$ and Desfontaines the values for $19\leq N\leq 20$. In this 
section, we mostly consider the functions $g(N,k)$ and 
$h(N,k)$ describing the general landscape of the pouring problem for $k\geq3$. 
In the following theorem, we present two monotony results.

\begin{theorem}
\label{TheMonotonicity}
For all positive integers $N$ and $k\geq3$, it holds:
\begin{itemize}
\item[\bf (a)] $g(N,k)\leq g(N+1,k)$, 
\item[\bf (b)] $\frac{k+1}{k} g(N,k)\leq g(N,k+1)$. 
\end{itemize}
\end{theorem}
\begin{proof}
\begin{itemize}
\item[\bf (a)]
Let $A=(a_1,a_2,\dots, a_k)$ be an initial 
state which is not $N$-pourable and $\sum_{i=1}^ka_i=g(N+1,k)$. By the 
definition of $N$-pourable, 
the state $A$ is not $(N-1)$-pourable. Therefore, $g(N,k)\leq \sum_{i=1}^ka_i=g(N+1,k)$.

\item[\bf (b)] 
Suppose to the contrary that for some positive integers $N$ and $k\geq3$, 
we have $g(N,k+1)< \frac{k+1}{k} g(N,k)$. Let $A=(a_1,a_2,\dots, a_{k+1})$ be an 
initial state which is not $(N-1)$-pourable, where $a_i\leq a_{i+1}$ for $ 1 \le i  \le k $, and 
$\sum_{i=1}^{k+1}a_i=g(N,k+1)$. In particular, we have $(k+1) a_{k+1}\geq 
g(N,k+1)$ and thus
\begin{eqnarray}
\label{ak_lb}
a_{k+1}\geq \frac{g(N,k+1)}{k+1}. 
\end{eqnarray}

Denote $n'=g(N,k+1)-a_{k+1}$. 
Consider the first $k$ vessels of $A$ as the initial state $B$. Notice that $B$ is 
not $(N-1)$-pourable, since we can apply 
the same pourings also to the state $A$ with $k+1$ vessels.  
The following holds:
\begin{eqnarray}
\label{eqb}
g(N,k) 
\le n' = g(N,k+1) -a_{k+1} \le \frac{k}{k+1} \cdot g(N,k+1) 
< \frac{k}{k+1} \frac{k+1}{k} g(N,k)=g(N,k),
\end{eqnarray}
where the first inequality
follows from the definitions of $ g(N,k)$ and $n'$, the second equality 
follows from the definition of $n'$, the third inequality follows 
from Inequality~\eqref{ak_lb}.  
and the fourth inequality from assumption.
By~\eqref{eqb}, we have a contradiction. 
It follows that $\frac{k+1}{k} g(N,k)\leq g(N,k+1)$.\qedhere
\end{itemize}
\end{proof}

 In the following two theorems, we give the exact value $g(N,k) $ 
for $ N \in \{ 1,2\}$  and for all integers $k \ge 3 $.

\begin{theorem}\label{theN=1}
For all integers $k \ge 3 $, it holds that $g(1,k)=k$. 
\end{theorem}
\begin{proof}
Let $k \ge 3 $ be given.
Since we have $N=1$, each vessel is non-empty in the initial 
state. The smallest state fulfilling this, is
the state $ (1,1,\dots,1)$.
This leads to $g(1,k)= k$. 
\end{proof}

\begin{theorem}\label{theN=2}
For all integers $k \ge 3 $, it holds that
$g(2,k)= \frac{k(k+1)}{2}$. 
\end{theorem}
\begin{proof}
Let $k \ge 3 $ be given.
Since we have $N=2$, each vessel has a unique non-zero value in the initial 
state. The smallest state fulfilling this, is
the state $(1,2,\dots,k)$. 
This leads to $g(2,k)= \frac{k(k+1)}{2}$. 
\end{proof}

We continue by giving an upper bound for $g(3,k)$. Note that by 
Theorems~\ref{TheMonotonicity}(a) and \ref{theN=2}, we have $g(3,k)\geq 
g(2,k)=\frac{k(k+1)}{2}$. Hence, the following theorem implies that 
$\frac{k(k+1)}{2}\leq g(3,k)\leq \lfloor \frac{5}{4}k^2 \rfloor $.  Furthermore, 
by Table~\ref{tab:gnk}, the upper bound of the following theorem is tight when 
$3\leq k\leq 8$.

\begin{theorem}
\label{propN=3Upp}
For all integers $k \ge 3 $, it holds that
$g(3,k)\leq \left\lfloor \frac{5}{4}k^2 \right\rfloor $.
\end{theorem}
\begin{proof}
We prove the assertion by giving an exactly $3$-pourable initial state with 
$\frac{5}{4}k^2 - \frac{1}{4} $ as the total sum of the vessels, for odd $k \ge 3$, and $\frac{5}{4}k^2$ as the total sum of the vessels, 
for even $k \ge 4$. 
Consider a state $A=(a_1,a_2,\dots,a_k)$ where $a_1=1$, $a_2=4$, 
$a_i=a_{i-2}+5$ for $i\geq3$.
First, we prove the following claims by induction:

\begin{description}

\item[\bf Claim 1:] For $ i \in \N $ it holds that
$ a_i = 
\left\{
\begin{array}{ll}
\frac{5}{2} i - \frac{3}{2} & \mbox{for odd $i$,} 
\\[0.6em]
\frac{5}{2} i - 1 & \mbox{for even $i$.} 
\end{array}
\right. $

\item[\bf Proof (Claim 1):] As induction base we consider the cases $i=1$ and $i=2$. These cases are 
clear, as $a_1=1$ and $a_2=4$.
The induction step is clear as well, as the $i$-th term and the 
$(i+2)$-th term in state $A$ differ exactly by $5$.

\item[\bf Claim 2:] It holds that
$ \sum_{s=1}^k a_s = 
\left\{
\begin{array}{ll}
\frac{5}{4} k^2  - \frac{1}{4} & \mbox{for odd $k$,} 
\\[0.6em]
\frac{5}{4} k^2  & \mbox{for even $k$.} 
\end{array}
\right. $

\item[\bf Proof (Claim 2):] As induction base we consider the cases $k=1$ and $k=2$. These cases are 
clear, as $a_1=1$ and $a_1+a_2=5$.

For the induction step assume that Claim 2 holds up to $k-1$.
We show that it holds also for $k$.
It holds by Claim 1 for odd $k$: 
\begin{eqnarray*}
\sum_{s=1}^k a_s = 
\frac{5}{4}(k-1)^2 + \left(\frac{5}{2} k - \frac{3}{2}\right)
=\frac{5}{4}k^2 - \frac{5}{2}k + \frac{5}{4}  
+ \frac{5}{2} k - \frac{3}{2}
=\frac{5}{4}k^2 - \frac{1}{4}.
\end{eqnarray*}

It holds by Claim 1 for even $k$: 
\begin{eqnarray*}
\sum_{s=1}^k a_s = 
\frac{5}{4} (k-1)^2 - \frac{1}{4}
+ \left(\frac{5}{2} k - 1\right)
=\frac{5}{4}k^2 -\frac{5}{2}k + \frac{5}{4} - \frac{1}{4}+ 
\frac{5}{2}k - 1
=\frac{5}{4}k^2. 
\end{eqnarray*}

This shows Claim~2.

\end{description}
\medskip

Notice that for any integer $k\geq3$, we have $a_1=1, a_2=4, a_3=6$. We 
may pour twice from $a_3$ to $a_1$ (leading to the state $a_1=4$, $a_2=4$, $a_3=3$) 
and then once from $a_2$ to $a_1$ to obtain an empty vessel. Hence, the state $A$ is 
$3$-pourable. 

Thus, it is left to show that at least $ N=3 $  pourings are required to empty a vessel. 
We show this based on the observation that after any single pouring from the initial state, 
no two vessels have an equal value 
and thus, after the first pouring, we still require at least two more pourings. 

Clearly, in the initial state all vessels have values $1,4,6,9 \! \mod 10$.
We distinguish between the following four cases.
\begin{itemize}
\item We pour from a vessel $a'$ into a vessel $a$ with value $ 1 \! \mod 10 $.

After the pouring, vessel $a$ has value $ 2 \! \mod 10 $.
Vessel $a'$ may have values $ 0,3,5,8 \! \mod 10 $.
None of these possible modulo values is equal to the remaining values $ 1,2,4,6,9 \! \mod 10 $.

\item We pour from a vessel $a'$ into a vessel $a$ with value $ 4 \! \mod 10 $.

After the pouring, vessel $a$ has value $ 8 \! \mod 10 $.
Vessel $a'$ may have values $ 7,0,2,5 \! \mod 10 $.
None of these possible modulo values is equal to the remaining values $ 1,4,6,8,9 \! \mod 10 $.

\item We pour from a vessel $a'$ into a vessel $a$ with value $ 6 \! \mod 10 $.

After the pouring, vessel $a$ has value $ 2 \! \mod 10 $.
Vessel $a$ may have values $ 5,8,0,3 \! \mod 10 $.
None of these possible modulo values is equal to the remaining values $ 1,2,4,6,9 \! \mod 10 $.

\item We pour from a vessel $a'$ into a vessel $a$ with value $ 9 \! \mod 10 $.

After the pouring, vessel $a$ has value $ 8 \! \mod 10 $.
Vessel $a$ may have values $ 2,5,7,0 \! \mod 10 $.
None of these possible modulo values is equal to the remaining values $ 1,4,6,8,9 \! \mod 10 $.

\end{itemize}
This shows the assertion.
\end{proof}

We use the following definition and the following lemma to give a lower bound for $g(4,k)$.

\begin{definition}
\label{def_seq}
Define the following two sequences $ (a_i)_{i \in \N} $ 
and $ (b_i)_{i \in \N} $ as follows.
\begin{itemize}
\item[{\bf (a)}]
Set $ a_1 := 1 $ and $ a_{i+1} := a_i + \lceil i/2 \rceil $ for $ i \ge 2 $,
i.e.,  $ (1,2,3,5,7,10,13,17,21,\dots) $.
\item[{\bf (b)}]
Set $ b_i := \sum_{j=1}^i a_j $ for $i \in \N  $,
i.e.,  $ (1,3,6,11,18,28,41,58,79,\dots) $.
\end{itemize}
\end{definition}

\begin{lemma}
\label{lemg4k}
\begin{itemize}
\item[{\bf (a)}]
For $ i \in \N $, it holds that $ a_i := \lfloor i^2/4 \rfloor + 1 $.
\item[{\bf (b)}]
For $ i \in \N $, it holds that $ b_i := 
\lfloor (1/12) i^3 + (1/8) i^2 + (11/12) i \rfloor $.
\end{itemize}
\end{lemma}

\begin{proof}
We prove the statements using induction on $ i \in \N $.
\begin{itemize}
\item[{\bf (a)}]
The cases $i=1$ and $i=2$ are clear,
as $ a_1 = 1 $ and $ a_2= 2 $ hold.
Let the assertion hold for $i \in \N $.
We will show that then it holds also for $ i+2$.
The claim holds by the following calculations:
\begin{eqnarray*}
a_{i+2} & = & a_i 
+ \lceil i/2 \rceil + \lceil (i+1)/2 \rceil 
\\
& = & 
\lfloor i^2/4 \rfloor + 1 + i + 1
\\
& = & 
\lfloor i^2/4 + i + 1 \rfloor + 1 
\\
& = & 
\lfloor (i+2)^2/4 \rfloor + 1.
\end{eqnarray*}
\item[{\bf (b)}]
The case $i=1$ 
is clear,
as $ b_1 = 1 $ holds.
Let the assertion hold for $i \in \N $.
We will show that then it holds also for $ i+1$.
It holds:
\begin{eqnarray}
\nonumber
b_{i+1} & = & b_i 
+ \lfloor (i+1)^2/4 \rfloor + 1
\\
& = & 
\label{i1_equal}
\lfloor (1/12) i^3 + (1/8) i^2 + (11/12) i \rfloor 
+ \lfloor (1/4) i^2 +(1/2)i + (1/4) \rfloor + 1
\\
\label{i2_equal}
& = & 
\lfloor (1/12) i^3 + (3/8) i^2 + (17/12) i + (9/8)  \rfloor 
\\
\nonumber
& = & 
\lfloor (1/12) i^3 + (1/4) i^2 + (1/4) i + (1/12)  
+ (1/8) i^2 + (1/4) i + (1/8)  
+ (11/12) i + (11/12)  
\rfloor 
\\
\nonumber
& = & 
\lfloor (1/12) (i+1)^3 + (1/8) (i+1)^2 + (11/12) (i+1) \rfloor. 
\end{eqnarray}
All equalities are clear except the one between \eqref{i1_equal} and \eqref{i2_equal}. We will show 
this equality in the following.
Set 
\begin{eqnarray*}
d_1 & := & (1/12) i^3 + (1/8) i^2 + (11/12) i,
\\
d_2 & := & (1/4) i^2 +(1/2)i + (1/4),
\\
d_3 & := & (1/12) i^3 + (3/8) i^2 + (17/12) i + (9/8).
\end{eqnarray*}
We further define for a real number $ r \in \R $, $ \alpha (r) = r - 
\lfloor r \rfloor \ge 0 $.

Note that it holds:
\begin{eqnarray}
\label{d123_rel}
d_3=d_1+d_2 +7/8. 
\end{eqnarray}

The equality of \eqref{i1_equal} 
and \eqref{i2_equal} is equivalent to:
\begin{eqnarray*}
\label{alpha_rel}
\lfloor d_1\rfloor+\lfloor d_2\rfloor+1 &=& \lfloor d_3\rfloor.
\end{eqnarray*}
This is equivalent to:
\begin{eqnarray*}
\label{alpha_rel2}
d_1+d_2+7/8 -(\alpha(d_1)+\alpha(d_2)) &=& d_3- (\alpha(d_3) +(1/8))
\end{eqnarray*}
and by \eqref{d123_rel} with
\begin{eqnarray}
\label{rel_alpha123}
\alpha(d_1) + \alpha(d_2) = \alpha (d_3) + (1/8).
\end{eqnarray}

Thus, it is left to show \eqref{rel_alpha123}.
For this, we distinguish two cases.
\begin{itemize}
\item $i$ is even.

Set $ i = 2a  $, where $ a \in \N $.

It holds:
\begin{eqnarray*}
d_1
& = & 
(1/12) (2a)^3 + (1/8) (2a)^2 + (11/12) (2a) 
\\
& = & 
(2/3) a^3 + (1/2) a^2 + (11/6) a
\\
& = & 
\frac{4 a^3 + 3 a^2 + 11 a}{6}.
\end{eqnarray*}

By checking all residue classes $ 0,1,2,3,4,5 \, \mod 6 $, we find that $ d_1 $ is an integer.
It follows that $ \alpha(d_1) = 0 $.

It holds:
\begin{eqnarray}
\nonumber
d_2
& = & 
(1/4) (2a)^2 + (1/2) (2a) + (1/4)
\\
\label{a_rel}
& = & 
a^2 + a + (1/4). 
\end{eqnarray}

It follows that $ \alpha (d_2) = 1/4 $.

By \eqref{d123_rel}
and \eqref{a_rel}, it follows:
\begin{eqnarray*}
d_3
& = & 
d_1 + d_2 + (7/8) 
\\
& = & 
d_1 + a^2 + a + (9/8). 
\end{eqnarray*}

It follows that $ \alpha (d_3) = 1/8 $.

\item $i$ is odd.

Set $ i = 2a  + 1 $, where $ a \in \N_0 $.

It holds:
\begin{eqnarray*}
d_1
& = & 
(1/12) (2a+1)^3 + (1/8) (2a+1)^2 + (11/12) (2a+1) 
\\
& = & 
(1/12) (8a^3 + 12a^2 + 6a +1 ) + (1/8) (4a^2 +4a +1) + (11/6) a + (11/12) 
\\
& = & 
(2/3) a^3 + a^2 + (1/2) a + (1/12) + (1/2) a^2 + (1/2) a + (1/8) + (11/6) a + (11/12) 
\\
& = & 
(2/3) a^3 + (3/2) a^2 + (17/6) a 
+ (9/8) 
\\
& = & 
\frac{4 a^3 + 9 a^2 + 17 a}{6} + (9/8).
\end{eqnarray*}

By checking all residue classes $ 0,1,2,3,4,5 \, \mod 6 $, we find that  
$ \frac{4 a^3 + 9 a^2 + 17 a}{6} $
is an integer.
It follows that $ \alpha(d_1) = 1/8 $.

It holds:
\begin{eqnarray*}
d_2
& = & 
(1/4) (2a+1)^2 + (1/2) (2a+1) + (1/4) 
\\
& = & 
a^2 + a + (1/4) + a + (1/2) + (1/4) 
\\
& = & 
a^2 + 2a+1.
\end{eqnarray*}

It follows that $ \alpha (d_2) = 0 $.

It holds:
\begin{eqnarray*}
d_3
& = & 
d_1 + d_2 + (7/8) 
\\
& = & 
(d_1-\alpha(d_1))+\alpha(d_1)  + d_2 + (7/8) 
\\
& = & 
(d_1-\alpha(d_1))+d_2+1.
\end{eqnarray*}

It follows that $ \alpha (d_3) = 0 $.

\end{itemize}

Thus, in both cases we have shown \eqref{alpha_rel}, and the proof is finished.\qedhere
\end{itemize}
\end{proof}

Now, we are ready to give a lower bound for $g(4,k)$.

\begin{theorem}
\label{prop=4Low}
It holds that $g(4,k)\geq \left\lfloor \frac{2k^3 +3k^2 + 22k}{24} \right\rfloor $.
\end{theorem}

\begin{proof}
Let $A=(a_1,a_2,\dots,a_k)$, with $a_i\leq a_{i+1}$ 
for $1\leq i\leq k-1$ and 
corresponding vessels $A_i$, be a state which is not $3$-pourable. In 
particular, this implies the following two conditions:
\begin{itemize}
\item $a_i\neq a_{i+1}$ for $1\leq i\leq k-1$, 
\item for each integer $j>0$, a gap $a_i-a_j$, for any $1\leq j<i\leq k$, can be used at most 
twice and then only if one vessel is involved in the gap twice, i.e., 
$a_{i_1}-a_{j_1} = a_{i_2}-a_{j_2}$, for any $1\leq i_1 < j_1 \le k$
and $1\le i_2 < j_2 \le k$ means automatically $i_1=j_2$ or $i_2=j_1$.
\end{itemize}
The first condition is clear. Assume that the second condition does not hold.
Then $a_{h}=a_i+j$ for some $j>0,h>i$ and $a_{p}= a_{q}+j$ for $p>q$ and $p>h$. 
Then we can pour from $A_{h}$ to $A_i$, and from  $A_{p}$ to 
$A_q$. At this point, the vessels $A_{h}$ and $A_{p}$ have equal value $j$, and 
finally we obtain an empty vessel by pouring from $A_{h}$ to $A_{p}$. 
Now we define a state $ (a_1,a_2,\dots, a_k) $, which fulfills both conditions and minimizes the value of each vessel. 
This is exactly the sequence $ (a_i)_{i \in \N} $ of Definition \ref{def_seq}(a).
It follows by Lemma~\ref{lemg4k}(b): 
\begin{eqnarray*}
g(4,k)
\; \geq \; b_k \; = \;  
\sum_{i=1}^k a_i
\; = \;
\left\lfloor \frac{2k^3 +3k^2 + 22k}{24} \right\rfloor. 
\end{eqnarray*}
\end{proof}

Note that together Theorems \ref{theN=1} and \ref{theN=2} and Theorems 
\ref{propN=3Upp} and \ref{prop=4Low} imply that $g(1,k)\in \Theta(k)$, $g(2,k)$, 
$g(3,k)\in \Theta(k^2)$ while $g(4,k)\in \Omega(k^3)$.

In Tables~\ref{tab:gnk} and \ref{tab:hnk}, we present some values and lower 
bounds for $g(N,k)$ and $h(N,k)$. 
In particular, we have computed for $k=3$ each case with $N\leq 65535$ 
(confirming the results of~\cite[A256001]{oeis25}), 
for $k=4$ each case with $N\leq 2047$,  
for $k=5$ each case with $N\leq 255$,  for $k=6$ each case with $N\leq 134$,  
for $k=7$ each case with $N\leq 100$ and  for $k=8$ each case with $N\leq 86$.
Note that by Corollary \ref{corhN3}, in Table~\ref{tab:hnk}, the values $h(N,3)$ are tight for 
$1\leq N\leq 9$ and that every other result in the table should be treated as only a 
lower bound although it seems highly likely that the values with small $N$ are tight. 

The values in the following tables are obtained by a C++ program and were run
on an x86 based parallel Linux cluster of the Christian-Albrechts University of 
Kiel, Germany. The total run was several core months.

\begin{table}[H]
\centering
\setlength{\tabcolsep}{0.5cm}
\begin{tabular}{l||c||c|c|c|c|c}
$N$ /	$k$	&	3	&	4 & 5 & 6 & 7 & 8	
\\ \hline
1	&	\phantom{00}3	&	\phantom{000}4 & \phantom{00}5 & \phantom{00}6 & \phantom{0}7 & \phantom{0}8	
\\ 
2	&	\phantom{00}6 & \phantom{00}10 & \phantom{0}15 & \phantom{0}21 & 28 & 36	
\\ 
3	&	\phantom{0}11	&	\phantom{00}20 & \phantom{0}31 & \phantom{0}45 & 61 & 80	
\\ 
4	&	\phantom{0}15	&	\phantom{00}40 & \phantom{0}71 & 123 & & 
\\ 
5	&	\phantom{0}23	&	\phantom{00}76 & 176 & & & 
\\ 
6	&	\phantom{0}27 & \phantom{00}177 & & & & 
\\ 
7	&	\phantom{0}45 & \phantom{0}387 & & & & 
\\ 
8	&\phantom{0}81 & \phantom{0}829	& & & & 
\\ 
9	&	105 & 1749	& & & & 
\end{tabular}
\caption{
Exact values $g(N,k) $ for some positive integers $N$ and $k \ge 3 $. 
\label{tab:gnk}}
\end{table}

\begin{table}[H]
\centering
\setlength{\tabcolsep}{0.5cm}
\begin{tabular}{l||c||c|c|c|c|c}
$N$ /	$k$	&	3	&	4 & 5 & 6 & 7 & 8	
\\ \hline
1	&	\phantom{000}5	&	\phantom{00}9 & \phantom{0}14 & \phantom{0}20 & \phantom{0}27 & \phantom{0}35
\\ 
2	&	\phantom{00}10 & \phantom{0}19 & \phantom{0}30 & \phantom{0}48 & \phantom{0}60 & \phantom{0}83	
\\ 
3	&	\phantom{00}20	&	\phantom{0}39 & \phantom{0}72 & 134 &  & 	
\\ 
4	&	\phantom{00}40	&	\phantom{0}86 & 201 &  & & 
\\ 
5	&	\phantom{00}80	&	204 &  & & & 
\\ 
6	&	\phantom{0}160 &  & & & & 
\\ 
7	&	\phantom{0}320 &  & & & & 
\\ 
8	&	\phantom{0}640 &	& & & & 
\\
9	&	1280 &	& & & & 
\end{tabular}
\caption{
Values $h(N,k) $ for some positive integers $N$ and $k \ge 3 $ (exact values for 
$k=3$ and lower bounds for $k \ge 4$). 
\label{tab:hnk}}
\end{table}

We continue by giving a monotonicity result for $h(N,k)$.

\begin{theorem}
\label{ThehMonotonicity}
For all positive integers $N$ and $k\geq3$, it holds:
$h(N,k)\leq h(N+1,k)$, 
\end{theorem}
\begin{proof}
Let $n=h(N,k)$. By definition, each state $A=(a_1,a_2,\dots,a_k)$ is 
$N$-pourable. Hence, each state $A$ is also $(N+1)$-pourable. Thus, we have 
$h(N,k)=n\leq h(N+1,k)$.
\end{proof}

While we have managed to give some initial results for the functions $g$ and $h$, 
there is still a lot to know. This is in particular true for the variants $g'$ and 
$h'$ of the functions $g$ and $h$ presented in Definition \ref{defgh}. In 
particular, we do not know whether $g'$ and $h'$ always exist or whether 
$g'(N,k)=g(N,k)$ and $h'(N,k)=h(N,k)$ holds. 

\begin{Question}
\label{question}
Do the following properties hold
for each positive integers $N$ and $k\geq3$?

\begin{enumerate}
\item $g'(N,k)$ and $h'(N,k)$ exist/are finite. 
\item $g'(N,k)\leq g'(N+1,k)$ and $h'(N,k)\leq h'(N+1,k)$.
\item $g(N,k)=g'(N,k)$.
\item $h(N,k)=h'(N,k)$.
\item $h(N,k)\leq h(N,k+1)$. 
\item $g'(N,k)\leq g'(N,k+1)$ and $h'(N,k)\leq h'(N,k+1)$.
\end{enumerate}
\end{Question}


Observe that we can say something about the relationships between the four 
parameters, when $g'$ and $h'$ exist and are finite. 
\begin{theorem}
\label{propgg'h'hRelation}
If the positive integers $N$ and $k\geq3$ are such that the values $g'(N,k)$ and $h'(N,k)$ exist and 
are finite, then it holds: 
$$g(N,k)\leq g'(N,k)\leq h'(N,k)\leq h(N,k).$$
\end{theorem}
\begin{proof}
We prove each of the three inequalities one by one, from left to right. 

First of all, we have $g(N,k)\leq g'(N,k)$, 
as the definition of $g'$ contains one more condition than the one of $g$. 

Secondly, we have $g'(N,k)\leq h'(N,k)$ since if the state $A=(a_1,a_2,\dots,a_k)$ 
with $a_1+a_2+\cdots+a_k=n$ is exactly $N$-pourable and each state with total 
value $n$ is $N$-pourable, then $g'(N,k)\leq n \le h'(N,k) $.

Finally, we have $h'(N,k)\leq h(N,k)$,
since the definition of $h'$ contains one more condition than the one of $h$. 
\end{proof}

\section{Lower bound for at least four vessels}
\label{sec_lower}

In Theorem \ref{thek3hBound}, we have seen that, in the case of three vessels, 
solving the pouring problem needs at least $\Omega(\log n)$ pourings. In this section, we generalize 
that result for the case of at least four vessels.

\begin{theorem}\label{thmLowerbound}
Let $t\geq3$ be an integer, 
and let $n$ be a sufficiently large integer. 
Then there exists a state $ (a_1,a_2,\dots,a_t) $ 
which is pourable in $\Omega(\log n)$ steps.
\end{theorem}
\begin{proof}
Consider a state $\left(a_1,a_2,\dots,a_t\right)=\left(n-\sum_{i=1}^{t-1} \left \lceil
n^{1/2^i} \right\rceil, \left\lceil n^{1/2^1}\right\rceil,\left\lceil 
n^{1/2^2} \right\rceil,\dots, \left\lceil
n^{1/2^{t-1}}\right\rceil \right)$ with respective vessels $A_1,A_2,\dots, A_t$. 
Notice that $\sum_{i=1}^ta_i=n$. Furthermore, if $n$
is sufficiently large, then we have $n\geq(t-1)\left\lceil \sqrt{n}  \, \right\rceil
>\sum_{i=1}^{t-1}\left\lceil n^{1/2^i}\right\rceil$ since
$t$ is a constant. Hence, the state contains only positive integers.

To make one of the vessels empty, we need to first  have two vessels with an equal value.
Consider the number of pourings we need for creating two 
vessels with equal value. Notice that if $n$ is sufficiently large, then we have: 
\begin{eqnarray}
\label{rel1}
\left\lceil n^{1/2^i}\right\rceil
>2\left\lceil n^{1/2^{i+1}}\right\rceil\ 
\mbox{for}\ i \in \{1,2,\dots, t-2\},
\\
\label{rel2}
n-\sum_{i=1}^{t-1}\left\lceil n^{1/2^i}\right\rceil
> 2 \left\lceil n^{1/2^1}\right\rceil.
\end{eqnarray}
Moreover, since  pouring from some vessel $X$ to another vessel $Y$, where
$|X|\geq2|Y|$, decreases the
value of $X$ to at least $|X|/2$ and at most doubles the value of vessel $Y$, we only need to
consider how many pouring steps we need to reach the state where the vessels $A_i$ and $A_{i+1}$
have  the same value for some $i\in \{1,2,\dots,t-1\} $, as other pairs of vessels
need even more pouring steps.

Let us first consider the case where $i\geq2$.
Suppose that we need exactly $j$ pourings for reaching the state, where the vessel $A_{i+1}$ has
value at least half as large as $A_i$. 
By Equality.~\eqref{rel1} it follows:
\begin{equation}\label{eqLower2j}
2^j\lceil n^{1/2^{i+1}} \rceil\geq \lceil n^{1/2^{i}}\rceil /2^{j+1}
\end{equation}
for $ i \in \{1,2,\dots,t-1\} $.
Inequality~(\ref{eqLower2j})  implies that $2^{2j+2} n^{1/2^{i+1}} \geq
n^{1/2^{i}}$ which is equivalent to $2j+2 \geq
\frac{1}{2^{i+1}}\log n $ and equivalent to $  j \geq \frac{1}{2^{i+2}}\log n-1$. Since $i\leq
t$, where $t$ is a constant, we need $j\in \Omega(\log n)$ pourings in this
case.

Let us next assume that we need exactly $j$ pourings for reaching the state,
where the vessel
$A_{2}$ has value at least half as large as $A_1$. 
By Equality.~\eqref{rel2} it follows:
\begin{equation}
\label{eqLower2ja12}
2^j\lceil n^{1/2^{1}} \rceil\geq \left(n-\sum_{i=1}^{t-1}\lceil n^{1/2^i}\rceil\right)/2^{j+1}.
\end{equation}
Inequality~(\ref{eqLower2ja12}) implies that $2^{2j+2} n^{1/2} \geq (n-\sum_{i=1}^{t-1}\lceil n^{1/2^i}\rceil) \geq n/2$ if
$n$ is sufficiently large. Hence, $j\geq \frac{1}{4}\log n-3/2$. Thus, we need $j\in \Omega(\log n)$ pourings also in this case.
\end{proof}

\section{Upper Bound for at least four vessels}
\label{sec_four}

In this section, we present Algorithm~\ref{alg4_vessel}
for solving the pouring problem for four vessels. 
This algorithm needs $\mathcal{O}(\log n\log\log n)$ pouring 
steps, improving on the best known algorithm for three vessels which requires $\mathcal{O}((\log n)^2)$ pouring 
steps (see Theorem~\ref{thmBoundPourings3}).
Note that a similar algorithm can be applied also for a greater number of vessels. 
	
\begin{algorithm}
	\caption{Four vessel pouring algorithm}
	\begin{algorithmic}[1]
	\Require  State $(a,b,c,d)$ with $a,b,c,d\in \N$.
	\Ensure New state $(a',b',c',d')$ in which one of $ a', b', c', d'$ is $0$.
\item Apply several times Algorithm~\ref{alg_frei} to the three vessels 
with smallest cardinality
\newline
until the smallest vessel has value smaller than  $\frac{n}{2\log n}$.
\label{lalg3_alg2a}

\item Order the vessels by their cardinality and call them $A$, $B$, $C$, $D$.
\label{lalg3_ord}

\Statex{(As after  this step, it is not any more poured into vessel $D$, but only out of $D$, $D$ is called \emph{pool} 
and \linebreak no longer \emph{vessel}.)}

\item Define $ e := e(A,B,C)$ such that 
$\gcd(|A|,|B|,|C|)$ is exactly $e$-even.
\label{lalg3_def}

\While{$\min\{|A|,|B|,|C|\}>0$ }
\label{lalg3_while_begin}

\If{$|A|\geq \frac{n}{4\log n}$}
\label{lalg3_if_largea}

\State Apply Algorithm~\ref{alg_frei} to the vessels $A$, $B$, $C$. 
\label{lalg3_alg2b}

\EndIf
\label{endifa}

\State Rename $A$, $B$ and $C$ so that $|A| \le  |B|\le |C|$. 
{\bf if} $|A| =  0$ {\bf then} {\bf Stop}.
\label{lalg3_rena}

\If{There are exactly two $(e+1)$-odd vessels}
\label{lalg3_if_odda}

\State Pour from one $(e+1)$-odd vessel to the other one.
\label{lalg3_pour_odd}

\Statex{\hspace*{2.5em} (If both vessels have equal value, pour always from $C$ to $B$, from $C$ to $A$, or from $B$ to $A$.)} 

\State Increase $e$ by $1$. 
\label{lalg3_inca}

\State Rename $A$, $B$ and $C$ so that $|A| \le  |B|\le |C|$. 
{\bf if} $|A| =  0$ {\bf then} {\bf Stop}.
\label{lalg3_renb}

\Else{\label{lalg3_else}
\If{Exactly three of the vessels are $(e+1)$-odd}
\label{lalg3_if_oddb}

\State Pour from $C$ to $B$.
\label{lalg3_pour_cb_one}

\State Rename $A$, $B$ and $C$ so that $|A| \le  |B|\le |C|$. 
{\bf if} $|A| =  0$ {\bf then} {\bf Stop}.
\label{lalg3_eif_oddb}

\EndIf
\label{endifb}

\If{$|C|$ is $(e+1)$-odd}
\label{lalg3_if_oddc}

\State Pour from $C$ to $B$ until $|C|<|B|$. 
\label{lalg3_pour_cb_sev}

\State Rename $A$, $B$ and $C$ so that $|A| \le  |B|\le |C|$. 
{\bf if} $|A| =  0$ {\bf then} {\bf Stop}.
\label{lalg3_rend}

\EndIf
\label{endifc}

\If{$|B|$ is $(e+1)$-odd and $|B| < \frac{n}{4\log n}$}
\label{lalg3_if_oddd}

\State Pour from $D$ to $B$.
\label{lalg3_pour_db}

\State Increase $e$ by $1$. 
\label{lalg3_incb}

\State Rename $A$, $B$ and $C$ so that $|A| \le  |B|\le |C|$. 
{\bf if} $|A| =  0$ {\bf then} {\bf Stop}.
\label{lalg3_rene}

\Else{
\If{$|B|$ is $(e+1)$-odd and $|B| \ge \frac{n}{4\log n}$}
\label{lalg3_if_odde}
\State Set $t_a=
\log\frac{\frac{n}{4\log n}}{|A|}$. 
\label{lalg3_set_ta}

\State Set $t_c$ such that $|C|$ is exactly $(e+t_c)$-even.
\label{lalg3_set_tc}

\State Set $t=\max \left\{ \min\{\left\lceil\frac{t_a}{2}\right\rceil,t_c\}, 1 \right\}$.
\label{lalg3_set_t}

\State Pour from $B$ to $A$ exactly $t-1$ times.
\label{lalg3_pour_ba}

\State Apply Algorithm~\ref{alg_janson} to the vessels $A$, $B$, $C$.
\label{lalg3_alg1}

\State Rename $A$, $B$ and $C$ so that $|A| \le  |B|\le |C|$. 
{\bf if} $|A| =  0$ {\bf then} {\bf Stop}.
\label{lalg3_renf}

\State Pour from $D$ to $A$ exactly $t$ times.
\label{lalg3_pour_da_sev}

\State Increase $e$ by $t$. 
\label{lalg3_incc}

\State Rename $A$, $B$ and $C$ so that $|A| \le  |B|\le |C|$. 
{\bf if} $|A| =  0$ {\bf then} {\bf Stop}.
\label{lalg3_reng}

\Else{
\If{$|A|$ is $(e+1)$-odd}
\label{lalg3_if_oddf}

\State Pour from $D$ to $A$.
\label{lalg3_pour_da_one}

\State Increase $e$ by $1$. 
\label{lalg3_incd}
\State Rename $A$, $B$ and $C$ so that $|A| \le  |B|\le |C|$. 
{\bf if} $|A| =  0$ {\bf then} {\bf Stop}.
\label{lalg3_renh}
\EndIf
}
\EndIf
}
\label{lalg3_endif}
\EndIf
}\EndIf
\EndWhile
\label{lalg3_while_end}
\end{algorithmic}
\label{alg4_vessel}
\end{algorithm}

\begin{theorem}\label{thm4vesselloglogAlg}
Let us have state $(a,b,c,d)$ and $a+b+c+d=n$. Then Algorithm~\ref{alg4_vessel} 
empties one of the vessels in $\mathcal{O}(\log n\log\log n)$ pourings.
\end{theorem}
\begin{proof}
In the following we sometimes write $A_i$/$B_i$/$C_i$ for the vessel $A$/$B$/$C$ 
before Step~$i$. If it is clear from the context, then we omit this 
index.
Let $ r:= \lceil \log n-2-\log \log n \rceil $ 
and $ r':= \lceil \log n \rceil $. 

Let us enumerate the while-loops by $w_i$, starting from $1$, where loop $w_i$
with index $i$ is the first while-loop at the beginning of which we have $e\geq i$. 
For example, $w_0=1$.
Observe in particular that we may have $w_{i+1}=w_{i+2}$ if we increase $e$
by more than one during a single while-loop (this may occur within the
if-clause of Step \ref{lalg3_if_odde} or at the beginning of the algorithm).

We divide the proof into two parts, namely into the correctness proof which shows 
that one of the vessels is emptied, and into the complexity proof that Algorithm 
\ref{alg4_vessel} needs 
$\mathcal{O}(\log n\log\log n)$ pourings.

\medskip

\noindent
{\bf Correctness:}
We show the following claims:

\begin{description}

\item[\bf (a)] After an increase of $e$ (in Steps~$\ref{lalg3_inca}$, 
\ref{lalg3_incb},~\ref{lalg3_incc},~\ref{lalg3_incd}), the 
value $\gcd(|A|,|B|,|C|)$ is exactly $e$-even.

\item[\bf (b)] In each round of the while-loop $e$ increases.

\item[\bf (c)] Steps~\ref{lalg3_pour_odd},~\ref{lalg3_pour_cb_one}, 
\ref{lalg3_pour_cb_sev},~\ref{lalg3_pour_ba} have valid pourings.

\item[\bf (d)] If we have done only valid pourings until Step $x$, 
then after Step $x$ the cardinality of the smallest vessel is smaller than 
$\frac{n}{2\log n}$.


\item[\bf (e)] If the first $r$ pourings from $D$, i.e., during 
Steps~\ref{lalg3_pour_db},~\ref{lalg3_pour_da_sev},~\ref{lalg3_pour_da_one} 
are poured into a vessel with cardinality smaller than $ \frac{n}{4 \log n} $, 
then during and after these $r$ pourings it holds that $|D| > \frac{n}{4 \log n} $.   

\item[\bf (f)] The first $r$ pourings from $D$
are poured into a vessel with cardinality smaller than $ \frac{n}{4 \log n} $. 

\item[\bf (g)] Steps~\ref{lalg3_pour_db},~\ref{lalg3_pour_da_sev}, 
\ref{lalg3_pour_da_one} have valid pourings for the $r$ first pourings from $D$.

\item[\bf (h)] After at most $r'$ increases of $e$, a vessel becomes empty.
\end{description}

\noindent
{\bf Proofs:}

\begin{description}
\item[\bf (a)] 
At the beginning of the while-loop we start with the situation
that $|A|$, $|B|$, $|C|$ are $e$-even, but $s$ of them for $1 \le s \le 3 $
are $(e+1)$-odd.

If $s=2$ (if-clause of Step~\ref{lalg3_if_odda}), then after pouring from one 
$(e+1)$-odd vessel to another one in Step~\ref{lalg3_pour_odd}, we obtain 
three $(e+1)$-even vessels (but not all of them are 
$(e+2)$-even),
i.e., $e$ is increased by $1$ in Step~\ref{lalg3_inca}.

If $s=3$ (if-clause of Step~\ref{lalg3_if_oddb}), then after pouring from $C$ 
to $B$, we obtain one $(e+1)$-odd vessel, namely $A$ (which may get renamed in 
Step~\ref{lalg3_eif_oddb}).

Thus, in the cases $s=1$ and $s=3$, we have exactly one $(e+1)$-odd vessel after 
Step~\ref{endifb}.
If $C$ is the $(e+1)$-odd vessel, then we pour from it to $B$ in 
Step~\ref{lalg3_pour_cb_sev} until 
the roles of B and C change. 
After Step~\ref{endifc}, $C$ is $(e+1)$-even, one
of $A$ or $B$ is exactly $e$-even, i.e., $(e+1)$-odd, 
and the other one is $(e+1)$-even.

We make a division into three cases, namely Cases 1, 2 for which $|B|$ is $ (e+1)$-odd 
(if-clauses of Steps~\ref{lalg3_if_oddd},~\ref{lalg3_if_odde}) and Case 3 for which $|A|$ is $ (e+1)$-odd
(if-clause of Step~\ref{lalg3_if_oddf}).

If $|B|$ is $(e+1)$-odd and $|B|<\frac{n}{4\log n}$ (if-clause of 
Step~\ref{lalg3_if_oddd}), 
then we pour from $D$ to $B$ making $|B|$ exactly $(e+1)$-even. Thus, $e$ is increased by $1$ in 
Step~\ref{lalg3_incb}. 

If $|B|$ is $(e+1)$-odd and $|B| \ge \frac{n}{4\log n}$ (if-clause of 
Step~\ref{lalg3_if_odde} ), 
then we pour from $D$ to $A$ exactly $t$ times in Step~\ref{lalg3_pour_da_sev}, where $t \ge 1 $ 
holds.

Observe that after Step~\ref{lalg3_set_t}, $|C|$ is 
$(e+t)$-even. After Step~\ref{lalg3_if_odde}, $|A|$ is $(e+1)$-even and after $t-1$ 
pourings in Step~\ref{lalg3_pour_ba} it becomes $(e+t)$-even. 
Thus, after Step~\ref{lalg3_pour_ba}, $|A|$ and $|C|$ are $(e+t)$-even. 
Furthermore, in Algorithm~\ref{alg_janson} we only pour into~$A$. Hence, after applying Algorithm~\ref{alg_janson}, $|A|$ and $|C|$ remain $(e+t)$-even.
As Algorithm~\ref{alg_janson} gives $B$ the smallest value,
$|B_{\ref{lalg3_pour_da_sev}}|$ and $|C_{\ref{lalg3_pour_da_sev}}|$ 
are $(e+t)$-even and $|A_{\ref{lalg3_pour_da_sev}}|$ is exactly $e$-even.

Thus, after pouring from vessel $D$ to 
vessel $A$ exactly $t$ times in Step~$\ref{lalg3_pour_da_sev}$,
the cardinalities of all vessels are $(e+t)$-even.
Thus, $e$ is increased by $t$ in Step~\ref{lalg3_incc}.

If $|A|$ is $(e+1)$-odd (if-clause of Step~\ref{lalg3_if_oddf}), 
then we pour from $D$ to $A$ making 
$|A|$ exactly  $(e+1)$-even. 
Thus, $e$ is increased by $1$ in Step~\ref{lalg3_incd}.

\item[\bf (b)] 
In each round of the while-loop one of the Steps~\ref{lalg3_inca}, 
\ref{lalg3_incb},~\ref{lalg3_incc},~\ref{lalg3_incd} is entered.
Each of these steps corresponds to an increase of $e$ by at least $1$.

\item[\bf (c)] 
Recall that a pouring is valid, if it is poured from a vessel 
into another one with smaller or equal cardinality. We make a case distinction for all pourings:
\begin{description}
\item[{\bf Step~\ref{lalg3_pour_odd}:}]

We can choose to pour from a vessel into another vessel with smaller or equal cardinality.

\item[{\bf Step~\ref{lalg3_pour_cb_one}:}]

The pouring is valid because $|B| \le |C| $. 

\item[{\bf Step~\ref{lalg3_pour_cb_sev}:}]

The pouring is valid as long $|B| \le |C| $ holds and stops when it no longer holds. 

\item[{\bf Step~\ref{lalg3_pour_ba}:}]

If $t=1$, then the claim holds since we do no pourings. 

If $t\geq2$, we have $ t_a \ge  t $. 
Then after the $t-1$ pourings from $B$ to $A$ we have:
\begin{eqnarray}
\nonumber
|A_{\ref{lalg3_alg1}}| & = & 2^{t-1}|A_{\ref{lalg3_pour_ba}}| \; \le \; 2^{t_a-1}
|A_{\ref{lalg3_pour_ba}}|=\frac{n}{8\log n}
\; < \; \frac{n}{4\log n} 
-\frac{n}{16\log n} 
\\
\label{eqB_alg1}
& \le & |B_{\ref{lalg3_pour_ba}}|  
-\frac{n}{16\log n} \; \le \; |B_{\ref{lalg3_alg1}}|.
\end{eqnarray}

It is left to show the last inequality of~\eqref{eqB_alg1}.
As during the $t-1$ pourings, 
exactly $ (2^{t-1} - 1 ) |A_{\ref{lalg3_pour_ba}}| $ 
is poured from $B_{\ref{lalg3_pour_ba}}$ to $A_{\ref{lalg3_pour_ba}}$,
we have to show that 
\begin{eqnarray}
\label{eq16}
(2^{t-1} - 1 ) |A_{\ref{lalg3_pour_ba}}| & \le & \frac{n}{16\log n}.
\end{eqnarray}

Inequality~\eqref{eq16}
is equivalent to
\begin{eqnarray*}
(2^{t-1} - 1 ) 
\frac{n}{2^{t_a} \cdot 4 \log n}
& \le & 
\frac{n}{16\log n}
\end{eqnarray*}
and equivalent to
\begin{eqnarray}
\label{eq4}
2^t - 2 & \le & 2^{t_a-1}.
\end{eqnarray}
As $t_a\geq2$ and $ t \le \lceil t_a/2 \rceil $ hold, 
Inequality~\eqref{eq4} follows from 
\begin{eqnarray*}
2^t - 2 \le 2^{\lceil t_a/2 \rceil } - 2 & \le & 2^{t_a-1}.
\end{eqnarray*}
Hence, Inequality (\ref{eq16}) holds.

\end{description}

\item[\bf (d)] Clearly, the claim holds after Step~\ref{lalg3_def} and thus it holds when 
we enter the while-loop. 
Consider then the while-loop. 
If we execute Step~\ref{lalg3_alg2b}, we 
decrease the value of the smallest vessel to a value smaller than 
$n/(4 \log n)$. Observe that after this, we may increase $|A|$ in 
Steps~\ref{lalg3_pour_odd},~\ref{lalg3_pour_ba},~\ref{lalg3_pour_da_sev}, 
\ref{lalg3_pour_da_one}. Moreover, 
we may either execute 
Step~\ref{lalg3_pour_odd} or Step~\ref{lalg3_pour_da_one}, or Steps 
\ref{lalg3_pour_ba} and~\ref{lalg3_pour_da_sev} (or none of them). 

\begin{description}
\item[{\bf Steps~\ref{lalg3_pour_odd} and \ref{lalg3_pour_da_one}:}]
In these steps we at most double the cardinality of $A$ from a value smaller than 
$n/(4 \log n)$ to a value smaller than $n/(2 \log n)$. After that we start with the 
next round of the while-loop and we again halve it 
to a value smaller than $n/(4 \log n)$ in Step~\ref{lalg3_alg2b}. 

\item[{\bf Steps~\ref{lalg3_pour_ba} and~\ref{lalg3_pour_da_sev}:}]
Since Step~\ref{lalg3_pour_ba}  is the first time in the while-loop 
where the value of the smallest 
vessel is increased, we have $t_a>0$ in Step 
\ref{lalg3_set_ta} since we have $|A|< n/(4 \log n)$ after Step~\ref{endifa}. 

If $t\geq2$ in Step~\ref{lalg3_set_t}, then after pouring in Step 
\ref{lalg3_pour_ba} we have 
$|A_{\ref{lalg3_alg1}}|=2^{t-1}|A_{\ref{lalg3_pour_ba}}| \leq \frac{n}{8\log n} 
< \frac{n}{2\log n}$. 
The second-last inequality is due to Inequality (\ref{eqB_alg1}).

Furthermore, 
we have 
$|A_{\ref{lalg3_alg1}}|=2^{t-1}|A_{\ref{lalg3_pour_ba}}|=\frac{n}{2^{t_a+1-t}4 
\log n}$. Since Algorithm~\ref{alg_janson} decreases the value of the smallest 
vessel, we have $|B_{\ref{lalg3_renf}}|< \frac{n}{2^{t_a+1-t}4 \log n}$. Thus, 
after we pour $t$ times to $A$ in Step~\ref{lalg3_pour_da_sev}, we have 
\begin{equation}\label{EqA_afterpour_da_sev}
|A_{35}|<\frac{n}{2^{t_a+1-2t}4 \log n}\leq \frac{n}{2 \log n}
\end{equation} since $t\leq \lceil t_a/2\rceil$.

\end{description}

\item[\bf (e)] 

Let the first $r$ pourings from $D$ be poured into a vessel with cardinality 
smaller than $ \frac{n}{4 \log n} $. 
Observe that after Step~\ref{lalg3_ord} it holds that $|D| >\frac{n}{4} $.
Each valid pouring from $D$ decreases $|D|$ by less than 
$ \frac{n}{4 \log n} $.
Thus, after $r$ pourings from $D$, it holds:
\begin{eqnarray*}
|D| \; > \; \frac{n}{4}-(\log n-1-\log \log n )\frac{n}{4\log n}
\; = \; \frac{n+n\log\log n}{4\log n}
\; > \; \frac{n}{4\log n}.
\end{eqnarray*}

\item[\bf (f)] 
We make again a case distinction:
\begin{description}
\item[{\bf Step~\ref{lalg3_pour_db}:}]

It holds that $|B| < \frac{n}{4 \log n} $, as we are in the if-clause of Step~\ref{lalg3_if_oddd}. 

\item[{\bf Step~\ref{lalg3_pour_da_sev}:}]

After Step~\ref{lalg3_pour_ba} we have $|A_{\ref{lalg3_alg1}} |= 
\frac{n}{2^{t_a+1-t}4\log n}\leq \frac{n}{4 \log n}$ 
since $t\geq1$. 
Furthermore, Algorithm~\ref{alg_janson} gives $B$ a value that is smaller than this,
and after renaming, vessel $B$ becomes $A$.
Thus, after the first $t-1$ pourings in Step~\ref{lalg3_pour_da_sev}, from $D$ 
to $A$, we have $|A_{\ref{lalg3_incc}}|<\frac{n}{2^{t_a+2-2t}4\log 
n}<\frac{n}{4\log n}$, as $ t_a+2-2t > 0 $ holds, which is 
equivalent to $ t < t_a/2 + 1 $. 
Hence, also the last pouring 
in Step~\ref{lalg3_pour_da_sev} satisfies the claim.


\item[{\bf Step~\ref{lalg3_pour_da_one}:}]

It holds that $|A| < \frac{n}{4 \log n} $ after Step~\ref{endifa} 
and this does not change before Step~\ref{lalg3_pour_da_one}. 
\end{description}

\item[\bf (g)] 
This follows directly from (e) and (f).

\item[\bf (h)] 
It holds:
\begin{eqnarray*}
2^{r'} \; = \; 
2^{\lceil \log n \rceil} 
\; \ge \; n  
\end{eqnarray*}
Thus, after $r'$ increases of $e$, $2^{r'}$ divides $|A|$, $|B|$, $|C|$.
This is only possible if at least one vessel is empty.
\end{description}

The correctness follows from {\bf (b)} and {\bf (h)}.

\medskip

\noindent
{\bf Complexity:}
\label{Alg3complexity}

We aim to give an upper bound for the number of pourings.
The following steps of Algorithm~\ref{alg4_vessel} contain pourings:
\ref{lalg3_alg2a},~\ref{lalg3_alg2b},~\ref{lalg3_pour_odd},~\ref{lalg3_pour_cb_one}, 
\ref{lalg3_pour_cb_sev},~\ref{lalg3_pour_db},~\ref{lalg3_pour_ba}, 
\ref{lalg3_alg1}, 
\ref{lalg3_pour_da_sev},~\ref{lalg3_pour_da_one}.
In the following we analyze how many pourings they contain.

\begin{description}
\item[Steps~\ref{lalg3_pour_odd},~\ref{lalg3_pour_cb_one}, 
\ref{lalg3_pour_db} and~\ref{lalg3_pour_da_one}:] 
We have only $1 \in \mathcal{O}(1)$, pourings.
\item[Step~\ref{lalg3_alg2a}:] 
We apply at most $q=\lceil \log\log n \rceil$ times 
Algorithm~\ref{alg_frei} where 
$\frac{n}{4}/2^q \le \frac{n}{2\log n}$. By Lemma~\ref{comp_jan_frei} 
each round uses at most  
$ \log \left( \frac{|B|}{|A|} \right) +2 <
\log \left( \frac{n}{|A|} \right) +2 \le
\log \left( 2 \log n \right) +2 =
\mathcal{O}\left(\log \log n\right)$ 
pourings if $ |A| \ge \frac{n}{2 \log n} $, 
and $ 0 \in \mathcal{O}\left(\log \log n\right)$ otherwise. 
Hence, this step works in $\mathcal{O}((\log\log n)^2)$ pourings.
\item[Step~\ref{lalg3_alg2b}:] 
By (b) and (h) of the correctness part, we have at most 
$r' \in \mathcal{O}\left(\log n\right)$ rounds of the
while-loop.
Analogously to the complexity proof of Step~\ref{lalg3_alg2a}, 
Algorithm~\ref{alg_frei} uses at most 
$\mathcal{O}\left(\log \log n\right)$ pourings.
In total, Step~\ref{lalg3_alg2b} has at most
$\mathcal{O}\left(\log n \log \log n\right)$ pourings.

\item[\bf Steps~\ref{lalg3_pour_cb_sev},~\ref{lalg3_pour_ba},~\ref{lalg3_alg1} and~\ref{lalg3_pour_da_sev}:]
We will use the following observations in the proofs of our subclaims.
\end{description}

\begin{observation}
\label{obs_a}
Consider a single iteration of the while-loop of 
Algorithm~\ref{lalg3_alg2a} after Step~\ref{endifc}, 
where the if-clause of Step~\ref{lalg3_if_odde} is not entered,
then $|A|$ at most doubles and we pour at most once into $A$.
\end{observation}

\begin{description}
\item[\it Proof of Observation~\ref{obs_a}.]
After Step~\ref{endifc}, $|A|$ might be changed only once,
namely doubled in Step~\ref{lalg3_pour_da_one}.
\end{description}

\begin{observation}
\label{obs_c}
In Algorithm~\ref{lalg3_alg2a}, it is never poured into vessel $C$. 
\end{observation}

\begin{description}
\item[\it Proof of Observation~\ref{obs_c}.]
As for $|B|=|C|$, we always pour from $C$ to $B$ (and analogous for $A$),
the only possibility to pour into vessel $C$ would be to pour from vessel $D$ 
which never occurs in Algorithm~\ref{lalg3_alg2a}. \end{description}

\begin{observation}
\label{obs_5_27}
Consider a single iteration of the while-loop of 
Algorithm~\ref{lalg3_alg2a}, where the if-clauses of Step~\ref{lalg3_if_largea} and
of Step~\ref{lalg3_if_odde} are not entered. Then it holds during 
this while-loop:
\begin{description}
\item[{\bf (a)}] 
$|A|$ at most doubles. 
\item[{\bf (b)}] 
$|C|$ reaches at least half of its previous value. 
\item[{\bf (c)}] 
The value of an $(e+1)$-even vessel 
cannot decrease and its parity can decrease only due to renaming.
\end{description}
\end{observation}

\begin{description}
\item[\it Proof of Observation~\ref{obs_5_27}.]
Consider a single iteration of the while-loop of 
Algorithm~\ref{lalg3_alg2a}, where the if-clauses in Step~\ref{lalg3_if_largea} and
in Step~\ref{lalg3_if_odde} are not entered. 
\begin{description}
\item[{\bf (a)}] 
We may pour into vessel $A$ only in Steps~\ref{lalg3_pour_odd} and~\ref{lalg3_pour_da_one} 
and only one of these may occur during each single run of the while-loop.
\item[{\bf (b)}] 
If we pour several times between two vessels $X$ and $Y$, by Observation~\ref{obs_bas},
the maximum value after all pourings is at least $ \max\left\{ |X|, |Y| \right\}/2 $.
Thus, the statement holds, if we enter Step~\ref{lalg3_pour_odd},
or if we enter Step~\ref{lalg3_pour_cb_one},
or if we enter Step~\ref{lalg3_pour_cb_sev},
or if we enter both Step~\ref{lalg3_pour_cb_one} and~\ref{lalg3_pour_cb_sev}.
If additionally we enter Step~\ref{lalg3_pour_db} or Step~ \ref{lalg3_pour_da_one},
clearly, the statement still holds.
\item[{\bf (c)}] The assertion follows, as no pourings from an 
$(e+1)$-even vessel are done (note that $D$ is not called vessel, but pool), 
if the if-clauses in Step~\ref{lalg3_if_largea} and
in Step~\ref{lalg3_if_odde} are not entered. \end{description}
\end{description}

\begin{description}
\item[\bf Step~\ref{lalg3_pour_cb_sev}:] 
Let $i \in \N $ with $0\le i\le r' = \lceil \log n\rceil$. 

Denote by $f_i$ the number of times we pour in Step~\ref{lalg3_pour_cb_sev} 
when $e=i$, and set $f_i=0$, if the algorithm does not reach 
Step~\ref{lalg3_pour_cb_sev} when $\gcd(|A|,|B|,|C|)$ is exactly $i$-even, or if it
never arrives to a state where the value $\gcd(|A|,|B|,|C|)$ is exactly $i$-even.


We will prove the following claim:

\begin{description}
\item[\bf Claim (i)] 
Let $i$ with $ 0 \le i \le r'$ and $f_i\ge \lfloor \log\log n \rfloor +8$, say 
\begin{eqnarray}
\label{eq_wi}
f_i= \lfloor \log\log n \rfloor +3+\ell_i, 
\end{eqnarray}
where $ \ell_i \ge 5 $ is an integer. 
Then for each integer $s$ with $1\leq s\leq \lfloor \ell_i/2\rfloor-1$
it holds that $f_{i+s} \le 1 $.
\end{description}


Notice first that $i+\lfloor\ell_i/2\rfloor-1\leq r'$ in the claim. 
This holds as after pouring $f_i$ times to the vessel $B$ in Step \ref{lalg3_pour_cb_sev}, 
we turn an $(i+1)$-even vessel into $(i+f_i+1)$-even vessel. Thus, 
$(i+f_i+1)\leq \log n$. In particular, this implies that 
$i+\lfloor\ell_i/2\rfloor-1< i+\ell_i<i+f_i<\log n\leq r'$. Thus, each $f_j$ is 
defined in the claim.

Observe that for each $i$ from the claim we have another set
$\{f_i,f_{i+1},\dots,f_{i+\lfloor \ell_i/2\rfloor-1}\}$, and all these sets
do not intersect. By the claim, the average value for each of these sets 
is at most $\frac{\lfloor\log\log n\rfloor+3+\ell_i+\lfloor 
\ell_i/2\rfloor-1}{\lfloor \ell_i/2\rfloor}\leq \frac{3\ell_i/2 + \log\log 
n+2}{\ell_i/4}< 
4\log\log n+14$. These leads to an average value over all $f_i$ contained
in one of these sets of $ \mathcal{O}(\log \log n)$.  
Observe further that for all other elements $f_i$, i.e., those which are not contained
in any of these sets, it holds that $f_i\le \log\log n+7$. Thus,
the average value over all these $f_i$ is also $ \mathcal{O}(\log \log n)$.  
So we have an overall average value of $ \mathcal{O}(\log \log n)$.  
By (b) and (h) of the correctness part, we have at most 
$r' \in  \mathcal{O}\left(\log n\right)$ rounds of the while-loop
and $\sum_{j=0}^{r'} f_j\in \mathcal{O}(\log n\log \log n)$ as claimed.
Thus, for Step~\ref{lalg3_pour_cb_sev}, it is only left to prove Claim (i).

\begin{description}
\item[\bf Proof of Claim (i)] 
We call the moment after Step \ref{endifc}, i.e., after we poured $f_i$ times in
Step \ref{lalg3_pour_cb_sev} as \emph{time moment 0}.

Immediately after $f_i - 1$  pourings in Step~\ref{lalg3_pour_cb_sev} it holds:
\begin{eqnarray*}
2^{f_{i}-1}|B_{\ref{lalg3_pour_cb_sev}}|| \; < \; |C_{\ref{lalg3_pour_cb_sev}}| - \left( 2^{f_{i}-1}- 1\right) 
|B_{\ref{lalg3_pour_cb_sev}}| 
\; < \; |C_{\ref{lalg3_pour_cb_sev}}| 
\end{eqnarray*}

Note again that $ |B_{\ref{lalg3_pour_cb_sev}}|| $
and $ |C_{\ref{lalg3_pour_cb_sev}}|| $ are the values of the vessels
$B$ and $C$ before Step~\ref{lalg3_pour_cb_sev}.

It follows by Equation~\eqref{eq_wi}:
\begin{equation}
\label{eqA19}
|A_{\ref{lalg3_pour_cb_sev}}|\leq 
|B_{\ref{lalg3_pour_cb_sev}}|<\frac{|C_{\ref{lalg3_pour_cb_sev}}||}{2^{\lfloor \log\log 
n \rfloor +2+\ell_i}}<\frac{n}{2^{\ell_i-1} 
4\log n}.
\end{equation}


In the following, we present four subclaims. 

In the first subclaim (i.1), we show how Algorithm \ref{alg4_vessel} 
advances from time moment $0$ when it does not enter the if-clause of Step \ref{lalg3_if_odde}.

\item[\bf Subclaim (i.1)] 
Let $s$ be an integer with 
$0\leq s\leq \lfloor \ell_i/2\rfloor -2$ such that at the beginning of the
while-loop $w_{i+s}$ we have $i+s \le e \le i+\lfloor \ell_i/2\rfloor -2$. Assume that
we do not enter the if-clause 
of Step \ref{lalg3_if_odde} during any of the while-loops $w_{i+s'}$ for $0\leq s' \leq s$.

Then the following statements hold at the end of the while-loop $w_{i+s'}$ 
for $0\leq s' \leq s$: 
\begin{description}
\item[\bf (1)] $|A|<\frac{n}{2^{\ell_i-2-s'}4\log n}$. 
\item[\bf (2)] $|A|\le \frac{|C|}{2^{f_i-1-2s'}}$.
\item[\bf (3)] Exactly one of the following two cases holds:
\begin{description}
\item[(I)] $|C|$ is $(i+f_i+1)$-even.
\item[(II)] $|B|$ is $(i+f_i+1)$-even and $|C|$ is $(i+f_i+1)$-odd. 
\end{description}
\item[\bf (4)]  $|C|\leq 2|B|$ in Case (II).
\end{description}
\item[\bf Proof of Subclaim (i.1)]
We prove this subclaim with induction on $s'$. Consider first the induction base. 
\begin{description}
\item[$\mathbf{s'=0}$:] 
We prove all four statements separately. We note that at the beginning of the 
while-loop $w_i$, we have $e=i$ since we enter Step \ref{lalg3_if_oddc} when 
$e=i$ by the definition of $f_i$.
\begin{description}

\item[\bf (1)] 
By Inequality (\ref{eqA19}), we have
$|A_{19}|<\frac{n}{2^{\ell_i-1}4\log n}$ right before 
time moment $0$. By our assumption, we do not enter the if-clause of Step
\ref{lalg3_if_odde}. 
Furthermore, $|A|$ at most doubles after time moment $0$, by Observation~\ref{obs_a}.
Thus, at the end of the while-loop $w_i$, we have 
\begin{eqnarray*}
|A|<2\cdot\frac{n}{2^{\ell_i-1}4\log n}
=\frac{n}{2^{\ell_i-2-s'}4\log n}. 
\end{eqnarray*}
This shows Statement~(1).

\item[\bf (2)] 
It holds that  $|A_{20}| \leq |A_{19}|\leq |B_{19}| = |C_{20}|/2^{f_i}$.
Furthermore, $|A|$ at most doubles after time moment $0$, by Observation~\ref{obs_a}.
Thus, at the end of the while-loop $w_i$, we have: 
\begin{eqnarray}
\label{eqfiac}
|A| \leq 2\cdot |A_{20}| \leq |C_{20}|/2^{fi-1} \leq |C|/2^{fi-1} = |C|/2^{f_i-1-2s'}
\end{eqnarray}
This shows Statement~(2). 

\item[\bf (3)] 
First, (I) and (II) cannot hold
simultaneously since then $|C|$ would be $(i+f_i+1)$-even and $(i+f_i+1)$-odd.
Furthermore, after time moment $0$, $|C|$ is $(i+f_i+1)$-even. Since we do
not enter the if-clause of Step \ref{lalg3_if_odde}, we do not decrease the
parity of vessel $C$ by Observation~\ref{obs_5_27}(c). Hence, the only way in which neither
$|B|$ nor $|C|$ is $(i+f_i+1)$-even at the end of the while-loop is if both
vessels $A$ and $B$ become larger than $C$ and vessel $C$ is renamed to be
vessel~$A$. 
By Observation~\ref{obs_a}, Inequality~\eqref{eqfiac}
and as renaming only may decrease $|A|$, this is not possible.
Hence, $|B|$ or $|C|$ remains $(i+f_i+1)$-even and (3) follows.

\item[\bf (4)] 
$|C|$ is $(i+f_i+1)$-even directly after time moment $0$. 
By Observation~\ref{obs_5_27}(c), the parity of $C$ does not decrease 
after time moment $0$ in this while loop, as renaming is also not possible.
Hence, in Case (II), we have poured into $B$
after time moment $0$ during the while-loop  $w_i$ so that $|B|$ has
surpassed $|C|$. Since we do not enter the if-clause of Step
\ref{lalg3_if_odde}, we pour into $B$ at most once after time moment $0$, and this must happen
in Step \ref{lalg3_pour_db}.
Then at the end of the while-loop, we have $|C| \le 2|B|$ and (4) follows.
\end{description}

\item[$\mathbf{s'-1 \to s'}$:] 
Assume that the claim holds for each integer $s^\ast$ with $0\leq s^\ast \le 
s'-1 <s\leq \lfloor \ell_i/2\rfloor -2$. We may further assume that $s'$ is such 
that $w_{i+s'-1}$ and $w_{i+s'}$ are  different while-loops since if they are 
the same while-loop, then  the claim holds by induction assumption. 
Consider the while-loop $w_{i+s'-1}$ and assume that
we do not enter the if-clause of Step \ref{lalg3_if_odde} during this and previous
while-loops. By induction assumption, at
the end of the while-loop $w_{i+s'-1}$ each of the statements (1), (2), (3) and (4)
holds. Hence, they also hold at the beginning of the while-loop $w_{i+s'}$. 
Again, we prove all four statements separately.

\begin{description}
\item[\bf (1)] 
Statement (1) holding at the beginning of the while-loop $w_{i+s'}$ 
implies that we do not enter the if-clause of Step
\ref{lalg3_if_largea}. Since we do not enter the if-clause of Step
\ref{lalg3_if_odde} and by Observation~\ref{obs_5_27}(a), 
$|A|$ at most doubles. Let $a$ be $|A|$ at the 
beginning of the while-loop $w_{i+s'}$ 
and $a'$ be $|A|$ at the 
end of the while-loop $w_{i+s'}$.
It follows:
\begin{eqnarray*} 
a' \le 2a < \frac{2n}{2^{\ell_i-2-(s'-1)}4\log n} 
=  \frac{n}{2^{\ell_i-2-s'}4\log n}. 
\end{eqnarray*} 
Thus, Statement (1) holds at the end of the while-loop $w_{i+s'}$. 

\item[\bf (2)] 
Again, we do not enter the if-clause of Step~\ref{lalg3_if_largea}. 
Let $ a $ and $ c $ be $ |A| $ and $|C|$, respectively, at the beginning of the while-loop $w_{i+s'}$, 
and let $ a' $ and $ c' $ be $ |A| $ and $|C|$, respectively, at the end of the while-loop $w_{i+s'}$. 
By Observation~\ref{obs_5_27}(a), (b), it follows:
\begin{eqnarray*} 
a' \le 2a \le \frac{2c}{2^{f_i-1-2(s'-1)}}
\le  \frac{4c'}{4 \cdot 2^{f_i-1-2s'}}
=  \frac{c'}{2^{f_i-1-2s'}}.
\end{eqnarray*} 
Thus, Statement (2) holds at the end of the while-loop $w_{i+s'}$. 

\item[\bf (3)] 
Since $|C|$ cannot be simultaneously 
$(i+f_i+1)$-even and $(i+f_i+1)$-odd, it is sufficient to
prove that at the end of the while-loop $w_{i+s'}$ one of $|B|$ or
$|C|$ is $(i+f_i+1)$-even. By induction assumption, at the beginning of the while-loop
$w_{i+s'}$, one of these values is $(i+f_i+1)$-even. Since we do not
enter the if-clauses in Steps \ref{lalg3_if_largea} or \ref{lalg3_if_odde}, the
parity or value of an $(i+f_i+1)$-even vessel cannot decrease during the
while-loop $w_{i+s'}$ by Observation~\ref{obs_5_27}(c). Hence, one of the vessels $B$ or $C$ remains $(i+f_i+1)$-even
unless it is renamed as vessel $A$. This can happen only if 
$|A|$ surpasses the value of the $(i+f_i+1)$-even vessel.
As Statement~(2) holds at the beginning of the while-loop $w_{i+s'}$, it 
follows that  $|A| \le \frac{|C|}{2^{f_i-1-2(s'-1)}}$. 

If the $(i+f_i+1)$-even vessel at the beginning of the while-loop was $C$,
then $|A|$ cannot surpass $|C|$, as $|A|$ at most doubles during the while-loop by Observation~\ref{obs_5_27}(a). 
If the $(i+f_i+1)$-even vessel at the beginning of the while-loop was $B$,
$|A|$ cannot surpass $|B|$, as by Statement (4) we had $|C|\leq 2|B|$ and thus we have at the beginning of the while-loop:
\begin{eqnarray*}
|A| \le \frac{|C|}{2^{f_i-1-2(s'-1)}}\leq
\frac{2|B|}{2 \cdot 2^{f_i-2s'}} 
=\frac{|B|}{2^{f_i-2s'}}.
\end{eqnarray*}
As $A$ at most doubles during the considered while-loop,
Statement (3) follows.

\item[\bf (4)] 
Assume that at the end of the while-loop $w_{i+s'}$, we are in Case (II). 

First, assume that at the beginning of the while-loop
$w_{i+s'}$ we are in Case (I). Recall that by Observation~\ref{obs_5_27}(c),
neither the value nor parity of vessel $C$ decreases. Since we are in Case (II)
at the end of the while-loop, $|B|$ surpasses $|C|$
during the while-loop $w_{i+s'}$. As by the explanations from proof of Statement (3), we see that
$|A|$ cannot surpass $|C|$, and 
$|B|$ surpasses $|C|$ by pouring into $B$ in Step~\ref{lalg3_pour_db}. 
Thus, $|B|$ at most doubles during the
while-loop $w_{i+s'}$. Furthermore, by Observation~\ref{obs_c}, 
we never pour into vessel $C$ and hence, $|C|\leq 2|B|$ at the end of 
the while-loop $w_{i+s'}$. 

Second, assume that at the beginning of the while-loop $w_{i+s'}$ we are in Case (II). 
By induction assumption, we have $|C|\leq 2|B|$ at the beginning of the
while-loop $w_{i+s'}$. By Observation~\ref{obs_c}, we do not pour into vessel $C$.
Furthermore, by 
Observation~\ref{obs_5_27}(c), the value or parity of $B$ does not decrease when it is
$(i+f_i+1)$-even. Since we are in Case (II) at the end of the 
while-loop $w_{i+s'}$, we have two possibilities.
As first possibility, we do not enter the if-clause of Step~\ref{lalg3_if_oddc}. Hence, vessel $B$ remains
$(i+f_i+1)$-even throughout the while-loop $w_{i+s'}$. 
As second possibility, we enter the if-clause of Step~\ref{lalg3_if_oddc} and after that
pour from $D$ into $B$ in Step \ref{lalg3_pour_db}. 
As $|A|$ does not surpass $|C|$, in both cases $|C| \le 2|B| $ holds at the end
of the while-loop $w_{i+s'}$. 
Thus, (4) follows.
\end{description}
\end{description}

In the second subclaim (i.2) we conclude the case where we do not enter 
the if-clause of Step \ref{lalg3_if_odde}.

\item[\bf Subclaim (i.2)] 
Let $s$ be an integer with 
$0\leq s\leq \lfloor \ell_i/2\rfloor -2$. Assume that we do not enter the if-clause 
of Step \ref{lalg3_if_odde} during any of the while-loops $w_{i+s'}$ for $0\leq s' \leq s$.
Then it holds that $ f_{i+s'+1} \le 1 $ 
for each $0\leq s'\leq s$. 

\item[\bf Proof of Subclaim (i.2)] 
Consider some $s'$ with $0\leq s'\leq s$ such that at the beginning of the
while-loop $w_{i+s'}$, we have $i+s' \le e \le i+\lfloor \ell_i/2\rfloor -2$. Such an
$s'$ exists since at the beginning of the while-loop $w_i$, we have $e=i$ by
the definition of $f_i$. Since we do not enter the if-clause of Step
\ref{lalg3_if_odde} during any of the while-loops $w_{i+s''}$ for $s''\leq s'$,
the statements of Subclaim (i.1) hold at the end of the while-loop $w_{i+s'}$
and at the beginning of the while-loop $w_{i+s'+1}$. Note that if $e>i+s'+1$
at the beginning of the while-loop $w_{i+s'+1}$, then we have 
$w_{i+s'+1}=w_{i+s'+2}=\cdots= w_{e}$ and $f_{i+s'+1}=0$.
Thus, we may assume that at the beginning of the while-loop $w_{i+s'+1}$ we
have $e=i+s'+1$.
Consider the two cases (I) and (II). 

If Case (I) holds, then $|C|$ is $(i+f_i+1)$-even, and thus $ f_{i+s'+1} = 0 $
holds. 
If Case (II) holds, then $|C|$ is $(i+f_i+1)$-odd, and $|C| \le 2 |B| $ holds at
the beginning of 
the while-loops $w_{i+s'+1}$ for $0\leq s' \leq s$.
So we pour at most once in Step \ref{lalg3_pour_cb_sev} 
and thus $ f_{i+s'+1} = 1 $ holds.
%
%
\end{description}

In the third subclaim (i.3), we show that entering
the if-clause of Step \ref{lalg3_if_odde} when $e=i+s$, causes  (under certain conditions)
an increase of $e$ to at least $i+\lfloor \ell_i/2\rfloor$ (and hence,
$f_{i+s'}=0$ for each $s< s'< \lfloor \ell_i/2\rfloor$). 

\begin{description}
\item[\bf Subclaim (i.3)] 

Assume that in the while-loop $w_{i+s} $ for an integer $s$ 
with $0\leq s\leq \lfloor \ell_i/2\rfloor -2$, 
we enter the if-clause of Step~\ref{lalg3_if_odde}, when $e=i+s$,
$|C|$ is $(i+f_i+1)$-even and $|A|<\frac{n}{2^{\ell_i-1-s}4\log n}$. 
Then in Step \ref{lalg3_incc}, we update the value of $e$ to at least $i+\lfloor\ell_i/2\rfloor$.
 
\item[\bf Proof of Subclaim (i.3)]
From $|A|< \frac{n}{2^{\ell_i-s-1} 4\log n}$ and the definition of $t_a$, 
it follows that $t_a> \ell_i-s-1\geq1$. Furthermore, $e+t_c\geq i+f_i+1 =  i + 
\lfloor\log\log n\rfloor+\ell_i+4$. Since $e=i+s$, we have $t_c\geq \lfloor\log\log n\rfloor+\ell_i+4-s$. Thus, 
$t=\max \left\{ \min\left\{ \lceil \frac{t_a}{2} \rceil, t_c \right\}, \, 1\right\} \geq 
\max\left\{\lceil\frac{\ell_i-s-1}{2} \right\rceil,\, 1\}$. 
Therefore, in Step \ref{lalg3_incc}, we increase 
$e$ from $i+s$ by at least $\lceil\frac{\ell_i-s-1}{2} \rceil$. Observe that after this, the value of 
$e$ is at least $i+s+\lceil\frac{\ell_i-s-1}{2}\rceil=i+\lceil\frac{\ell_i+s-1}{2}\rceil\geq i+\lceil\frac{\ell_i-1}{2}\rceil= i+\lfloor\frac{\ell_i}{2}\rfloor$.
\end{description}

In the fourth subclaim (i.4), 
we show that if we enter the if-clause of Step \ref{lalg3_if_odde} for the first time when $e=i+s'$,
then we have $f_{i+s}\leq 1$ for each $s' \leq s \leq \lfloor \ell_i/2\rfloor -2$.

\begin{description}
\item[\bf Subclaim (i.4)] 
If we enter the if-clause of Step \ref{lalg3_if_odde} during the while-loop $w_{i+s'}$ for $0\leq s' \leq \lfloor \ell_i/2\rfloor -2$, 
then it holds that $f_{i+s+1}\le 1$ for each integer $s$ with $s'\leq s \leq \lfloor \ell_i/2\rfloor -2$.

\item[\bf Proof of Subclaim (i.4)] 

Assume that we enter the if-clause of Step \ref{lalg3_if_odde} during
the while-loop $w_{i+s'}$ the first time after time moment $0$. Note that if at the beginning of the while-loop $w_{i+s'}$ we have $i+s'<e$, 
then we have 
$f_{i+s'}=0$. Hence, we assume that $e=i+s'$.

First, consider the case with $s'=0$, i.e., we reach the if-clause of Step~\ref{lalg3_if_odde} in the same
while-loop as time moment $0$. By Inequality \eqref{eqA19}, right before
time moment $0$ we have $|A|< \frac{n}{2^{\ell_i-1}4\log n}$. This
is true also when we enter the if-clause of Step \ref{lalg3_if_odde}.
Furthermore, since we pour $f_i$ times in Step \ref{lalg3_pour_cb_sev} to an
$(i+1)$-even vessel, $|C|$ is $(i+f_i+1)$-even after time moment $0$-.
Hence, $C$ is also $(i+f_i+1)$-even  when we enter the if-clause of Step
\ref{lalg3_if_odde}. As all three conditions of (i.3) are fulfilled,, we can apply it for $s=0$.
By (i.3), we update the
value of $e$ to at least $i+\lfloor \ell_i/2\rfloor$ in Step \ref{lalg3_incc}.
Hence, $f_{i+s+1}=0$ for each $0\leq s\leq \lfloor \ell_i/2\rfloor -2$.

Second, consider the case where we enter the if-clause of Step \ref{lalg3_if_odde}
for some $0<s'\leq \lfloor \ell_i/2\rfloor -2$. All statements of 
(i.1) apply at the end of the while-loop $w_{i+s'-1}$ and hence, at the
beginning of the while-loop $w_{i+s'}$. 
Let us next show that the
conditions of (i.3) apply when we enter the if-clause of Step
\ref{lalg3_if_odde}. First of all, by (i.1)(1), we have
$|A|<\frac{n}{2^{\ell_i-1-s'}4\log n}$ at the beginning of the while-loop.
Hence, we do not enter the if-clause of Step \ref{lalg3_if_largea}. Furthermore,
as by (i.1)(3),  $|B|$ or $|C|$ is $(i+f_i+1)$-even, we do not enter the if-clause of Step
\ref{lalg3_if_oddb}. Thus, $|A|<\frac{n}{2^{\ell_i-1-s'}4\log n}$ when we enter
the if-clause of Step \ref{lalg3_if_odde}. If we are in Case (I) at the
beginning of the while-loop, then $|C|$ remains $(i+f_i+1)$-even until we enter
the if-clause of Step \ref{lalg3_if_odde}. If on the other hand we are in Case (II)
at the beginning of the while-loop, then we decrease the parity of $B$ at some
point between Steps \ref{lalg3_else} and \ref{lalg3_if_odde}. 
By Observation~\ref{obs_5_27}(c), this is possible only if we rename vessel $B$, 
i.e., if we enter
the if-clause of Step \ref{lalg3_if_oddc} which results to vessel $C$ becoming
$(i+f_i+1)$-even. 
Thus, the conditions of (i.3) hold and hence, we increase $e$ to at
least $i+\lfloor\ell_i/2\rfloor$ in Step \ref{lalg3_incc}. Hence, $f_{i+s}=0$ for
each $s'<s\leq \lfloor \ell_i/2\rfloor -1$. Therefore, (i.4) follows.
\end{description}

Together (i.2) and (i.4) show that each $f_{i+s}\leq 1$ for 
$1\leq s\leq \lfloor \ell_i/2\rfloor -1$ in both cases (if we 
enter the if-clause of Step \ref{lalg3_if_odde} 
for some integer $s'$ with $1\leq s'\leq \lfloor \ell_i/2\rfloor -1$ 
or if we do not enter it). Claim (i) follows.

\end{description}

Finally, we are left with the average number of pourings within the if-clause of 
Step~\ref{lalg3_if_odde}. There, we execute pourings in 
Steps~\ref{lalg3_pour_ba},~\ref{lalg3_alg1} and~\ref{lalg3_pour_da_sev}.

\begin{description}
\item[\bf Steps~\ref{lalg3_pour_ba} and~\ref{lalg3_pour_da_sev}:] 

First of all, observe that in Step \ref{lalg3_pour_ba} we pour $t-1$
times and in Step \ref{lalg3_pour_da_sev} we pour $t$ times and afterwards, in
Step \ref{lalg3_incc}, we increase $e$ by $t$. Hence, in both cases we pour at most once on
average for each increase of $e$. 

\item[\bf Steps~\ref{lalg3_alg1}:]

During a single while-loop 
by Lemma~\ref{comp_jan_frei}, 
at most $\log \left(\frac{|B|}{|A|}\right)+2$ pourings may occur. 
Since $|B|\leq |C|$, we have $|B|<n/2$. In Step
\ref{lalg3_set_ta}, we have $t_a=\log\frac{\frac{n}{4\log n}}{|A|}$ and after
that in Step \ref{lalg3_pour_ba}, $|A|$ is doubled $t-1$ times. Thus,
before Step \ref{lalg3_alg1}, we have
\begin{eqnarray*}
|A_{\ref{lalg3_alg1}}|= n/(2^{t_a-t+1}4\log n).
\end{eqnarray*} 
Therefore, we have at most 
\begin{eqnarray}
\nonumber
& & \log \left(\frac{|B_{\ref{lalg3_alg1}}|}{|A_{\ref{lalg3_alg1}}|}\right)+2\leq
\log \left( \frac{n/2}{n/(2^{t_a-t+1}4\log n)} \right) +2 
\\
\label{pour32}
& = & \log(4 \cdot 2^{t_a-t}\log 
n)+2=t_a-t+4+\log\log n
\end{eqnarray}
pourings.
Furthermore, we increase $e$ by $t$ in Step \ref{lalg3_incc}. Assume first
that $t_a \leq 10t +\log\log n$. In this case, for each increase of $e$, we
pour at most $\frac{t_a-t+4+\log\log n}{t}\leq \frac{9t+2\log\log n+4}{t}\leq
2\log\log n+13\in \mathcal{O}(\log\log n)$ times in Step \ref{lalg3_alg1}.
Hence, we assume from now on that 
\begin{eqnarray*}
t_a & > & 10t +\log\log n.
\end{eqnarray*} 
We denote by $t_{a,i}$ ($t_i$) the value of $t_a$ ($t$) during the while-loop $w_i$.



\begin{description}
\item[\bf Claim (j)] 
Let $i$ with $ 0 \le i \le r'$. Let us 
enter the if-clause of Step \ref{lalg3_if_odde} during the while-loop $w_i$ 
and $t_{a,i} > 10t_i +\log\log n
$, say 
\begin{eqnarray*}
t_{a,i} = 10t_i + a_i +\lfloor \log\log n\rfloor, 
\end{eqnarray*}
where $ a_i $ is a positive integer. 
Then for each integer $s$ with $1\leq s\leq t_i + t_{a,i}/4$
it holds that we do not enter the if-clause of Step \ref{lalg3_if_odde} during 
the while-loop $w_{i+s}$.\end{description}

Notice that if Claim (j) holds, then by Inequality~\eqref{pour32}, 
during the while-loops $w_{j}$ for $i\leq
j\leq i+t_i+t_{a,i}/4$ we pour in total $t_{a,i}-t_i+4+\log\log n$ times in
Step \ref{lalg3_alg1}. Hence, on average during these while-loops, 
there are at most
\begin{eqnarray*}
\frac{t_{a,i}-t_i+4+\log\log n}{t_i+t_{a,i}/4+1} \leq\frac{9t_i+2\log\log
n+a_i+4}{\frac{14t_i+\lfloor \log\log n\rfloor+a_i}{4}} \leq
4 \cdot \frac{9t_i+2\log\log n+a_i+4}{11t_i+\log\log n+a_i+2}< 8 
\end{eqnarray*}
pourings in Step~\ref{lalg3_alg1}. 

\begin{description}

\item[\bf Proof of Claim (j)] 

By Inequality~\eqref{EqA_afterpour_da_sev} and as we do not reach the 
if-clause of Step~\ref{lalg3_if_oddf}, it holds:
\begin{eqnarray}
\label{eqA36}
|A_{\ref{lalg3_endif}}| \le
|A_{\ref{lalg3_reng}}|=|A_{\ref{lalg3_incc}}|<\frac{n}{2^{t_{a,i}+1-2t_i}4\log n}.
\end{eqnarray}


Before Step \ref{lalg3_if_odde} it holds that $|C_{\ref{lalg3_if_odde}}|\geq
|B_{\ref{lalg3_if_odde}}|\geq \frac{n}{4\log n}$. Since
$|C_{\ref{lalg3_endif}}|+|B_{\ref{lalg3_endif}}|+|A_{\ref{lalg3_endif}}|>|C_{\ref{lalg3_if_odde}}|+|B_{\ref{lalg3_if_odde}}|+|A_{\ref{lalg3_if_odde}}|$,
it holds: 
\begin{eqnarray*}
|C_{\ref{lalg3_endif}}|+|B_{\ref{lalg3_endif}}|>
|C_{\ref{lalg3_if_odde}}|+|B_{\ref{lalg3_if_odde}}|+
|A_{\ref{lalg3_if_odde}}|-|A_{\ref{lalg3_endif}}| \ge
2 \cdot \frac{n}{4\log n}-\frac{n}{2^{t_{a,i}+1-2t_i}4\log n}. 
\end{eqnarray*}
The last inequality is due
to Inequality (\ref{eqA36}). 
It follows because of $t_{a,i}+2-2t_i > 0$:  
\begin{eqnarray}
\label{c8}
|C_{\ref{lalg3_endif}}|>\frac{n}{4\log
n}-\frac{n}{2^{t_{a,i}+2-2t_i}4\log n} \ge \frac{n}{8\log n}. 
\end{eqnarray}

As explained in the second last paragraph of Claim (a) of the proof of correctness,
$|A_{\ref{lalg3_pour_da_sev}}|$ is exactly $e$-even. 
Furthermore, $|A_{\ref{lalg3_reng}}|$ is exactly $e$-even since we
obtain $|A_{\ref{lalg3_reng}}|$ from 
$|A_{\ref{lalg3_pour_da_sev}}|$ by pouring exactly $t_i$ times and then 
increasing $e$ by $t_i$. 

Observe that since we increase $e$ by $t_i$ in Step
\ref{lalg3_incc}, the while-loops $w_{i+1},w_{i+2},\dots,w_{i+t_i}$ are the same.
Furthermore, at the beginning of the while-loop $w_{i+t_i}$, at least one of
$|A|$ or $|B|$ is exactly $e$-even and is smaller than
$$\frac{n}{2^{t_{a,i}+1-2t_i}4\log n}<\frac{n}{2^{t_{a,i}/2}4\log n}.$$ We fix
this vessel as vessel $X$. 

\item[\bf Subclaim (j.1)] 
For $1\leq s' \leq t_{a,i}/4+1$ the following holds: 
\begin{description}
\item[\bf (1)] $|X|<\frac{n}{2^{t_{a,i}+1-2t_i-s'}4\log n}$ 
at the end of while-loop $w_{i+t_i+s'-1}$. 
\item[\bf (2)] $|X|$ is exactly $e$-even
at the end of while-loop $w_{i+t_i+s'-1}$. 
\item[\bf (3)] $|C|>\frac{n}{2^{s'}8\log n}$ 
at the end of while-loop $w_{i+t_i+s'-1}$. 
\item[\bf (4)] We do not enter the if-clause of Step \ref{lalg3_if_odde} during
the while-loop $w_{i+t_i+s'-1}$.
\end{description}

\item[\bf Proof of Subclaim (j.1)] 
The proof is by induction on $s'$.
\begin{description}
\item[$\mathbf{s'=1}$:] 
In this case we consider the while-loop $w_{i+t_i}$. 
Since $$|X|<\frac{n}{2^{t_{a,i}+1-2t_i-(s'-1)}4\log n} 
<\frac{n}{2^{t_{a,i}/2}4\log n} < \frac{n}{4\log n},$$ we do not enter 
Step~\ref{lalg3_if_largea}. 
We consider three cases for the beginning of the while-loop $w_{i+t_i}$: 
\begin{itemize}
\item[-] Exactly three vessels are $(e+1)$-odd. 

We pour from $C$ to $B$ in Step \ref{lalg3_pour_cb_one} and
later once to vessel $A$ in Step \ref{lalg3_pour_da_one}. 

\item[-] Exactly two vessels are $(e+1)$-odd. 

One of these two vessels is $X$ and we pour exactly once into $X$ in Step~\ref{lalg3_pour_odd}. 

\item[-] Exactly one vessel is $(e+1)$-odd. 

The vessel being $(e+1)$-odd is $X$. If $X$ is
$A$, then we pour into $X$ in Step \ref{lalg3_pour_da_one} and if $X$ is $B$, then we
pour into $X$ in Step \ref{lalg3_pour_db}. 

\end{itemize}

In all three cases, we increase the parity of $X$ by one proving Statement (2) and we
pour into $X$ at most once proving Statement (1). Furthermore, we pour from $C$ at most once,
which by Observation~\ref{obs_bas}, after a possible renaming, proves Statement (3). 
Finally, we do not enter the if-clause of Step \ref{lalg3_if_odde} in any of these three cases,
proving Statement (4). 

Hence, (j.1) follows for $s'=1$.

\item[$\mathbf{s' \to s'+1}$:] 
Assume that (j.1) holds at the end of the while-loop
$w_{i+t_i+s'-1}$, where $1\leq s'\leq  t_{a,i}/4$. Hence, it
also holds at the beginning of the while-loop $w_{i+t_i+s'}$. Since
$|X|<\frac{n}{2^{t_{a,i}+1-2t_i-s'}4\log n}<\frac{n}{2^{t_{a,i}/2}4\log
n}<\frac{n}{4\log n}$, we do not enter Step \ref{lalg3_if_largea}. Notice that
at this point we still have $|X|<|C|$ since $|C|>\frac{n}{2^{s'}8\log
n}\geq\frac{n}{2^{t_{a,i}/4}8\log n}$ as $t_{a,i}/4< t_{a,i}/2 - 1$.
Thus, the situation is the same as in the induction base $s'=1$,
and the proof is analogous.
\end{description}

\end{description}
By Subclaim (j.1), Claim (j) follows.\qedhere
\end{description}
\end{proof}

\section{Conclusions and Future Work}
\label{sec_conc}

In this work, we have generalized the pouring problem from $k=3$ vessels to 
$k\geq2$ vessels. For $k=2$, we have shown that we may solve it
only for the initial states $(a,b)$ with $a+b=2^t$ for $t\in \N$. Furthermore, while 
in the case of $k=3$ vessels, the best known upper bound for the required number 
of pouring steps is in $\mathcal{O}((\log n)^2)$ (see 
Theorem~\ref{thmBoundPourings3}), we have managed to give an 
upper bound in $\mathcal{O}(\log n\log\log n)$ pouring steps for the case of at 
least four vessels.

However, it still remains an open problem whether either of these upper bounds 
can be improved. In particular, we only have a lower bound of $\Omega(\log n)$
for $ k \ge 3 $ vessels.
We have also introduced the notations $g(N,k)$ and $h(N,k)$ for easier 
discussion of related problems and better understanding of the landscape of 
pouring problems. While we have presented some initial results for these 
parameters, many open problems (see Question~\ref{question}) still remain for the two parameters. Another line of open questions is determining the asymptotic behaviours of $g(N,k)$ when $N$ is fixed and $k$ grows. We have show that for $g(1,k)\in \Theta(k)$, $g(2,k),g(3,k)\in \Theta(k^2)$ and $g(4,k)\in \Omega(k^3)$.

\section*{Acknowledgments}\label{sec_ackn}

\noindent
The research of Tuomo Lehtilä was partially funded by the Academy of Finland grants 338797 and 358718.

\noindent We would like to thank John Tromp for discussions and ideas to improve our C++ 
program for computing the values $ g(N,k) $ and $ h(N,k)$.

\bibliographystyle{plainnat}
\bibliography{main81}

\end{document}